\documentclass[12pt]{amsart}
\usepackage{amssymb,amsmath,color}
\usepackage{color}
\usepackage{tikz}
\usepackage{pgf}
\usepackage{bm}
\usepackage{tkz-graph}
\usetikzlibrary{positioning}
\usetikzlibrary{calc}
 \usepackage{array,multirow,graphicx}
 \usepackage{float}
 \usepackage{pdflscape}
 \usepackage[left=1in,right=1in,top=0.9in,bottom=0.9in]{geometry}
 \usepackage{url}
 \usepackage{pdfpages}

\newtheorem{theorem}{Theorem}[section]

\newtheorem{lemma}[theorem]{Lemma}
\newtheorem{corollary}[theorem]{Corollary}
\newtheorem{proposition}[theorem]{Proposition}
\theoremstyle{definition}
\newtheorem{definition}[theorem]{Definition}
\newtheorem{example}[theorem]{Example}

\theoremstyle{remark}

\usepackage{stackengine}
\newcommand\ledot{\mathrel{\ensurestackMath{%
  \stackengine{-.5ex}{\lessdot}{-}{U}{c}{F}{F}{S}}}}
  \newcommand\gedot{\mathrel{\ensurestackMath{%
  \stackengine{-.5ex}{\gtrdot}{-}{U}{c}{F}{F}{S}}}}

\newcommand{\Hzero}[3]{
\begin{tikzpicture}[baseline={([yshift=-.5ex]current bounding box.center)}]
\node[fill=black, circle, inner sep=1pt, minimum size=0.1cm,label=right:\footnotesize{#1}] (1) {};
\node[fill=black, circle, inner sep=1pt, minimum size=0.1cm, label=left:\footnotesize{#2}] (2) [below left = 0.3cm of 1] {};
\node[fill=black, circle, inner sep=1pt, minimum size=0.1cm, label=right:\footnotesize{#3}] (3) [below right = 0.3cm of 1] {};
\end{tikzpicture}}

\newcommand{\Hone}[3]{
\begin{tikzpicture}[baseline={([yshift=-.5ex]current bounding box.center)}]
\node[fill=black, circle, inner sep=1pt, minimum size=0.1cm, label=right:\footnotesize{#1}] (1) {};
\node[fill=black, circle, inner sep=1pt, minimum size=0.1cm, label=left:\footnotesize{#2}] (2) [below left = 0.3cm of 1] {};
\node[fill=black, circle, inner sep=1pt, minimum size=0.1cm, label=right:\footnotesize{#3}] (3) [below right = 0.3cm of 1] {};
\draw (2)--(1);
\end{tikzpicture}}

\newcommand{\Htwo}[3]{
\begin{tikzpicture}[baseline={([yshift=-.5ex]current bounding box.center)}]
\node[fill=black, circle, inner sep=1pt, minimum size=0.1cm, label=right:\footnotesize{#1}] (1) {};
\node[fill=black, circle, inner sep=1pt, minimum size=0.1cm, label=left:\footnotesize{#2}] (2) [below left = 0.3cm of 1] {};
\node[fill=black, circle, inner sep=1pt, minimum size=0.1cm, label=right:\footnotesize{#3}] (3) [below right =0.3cm of 1] {};
\draw (3)--(1)--(2);
\end{tikzpicture}}

\newcommand{\uHtwo}[0]{
\begin{tikzpicture}[baseline={([yshift=-.5ex]current bounding box.center)}]
\node[fill=black, circle, inner sep=1pt, minimum size=0.1cm] (1) {};
\node[fill=black, circle, inner sep=1pt, minimum size=0.1cm] (2) [below left = 0.3cm of 1] {};
\node[fill=black, circle, inner sep=1pt, minimum size=0.1cm] (3) [below right =0.3cm of 1] {};
\draw (3)--(1)--(2);
\end{tikzpicture}}

\newcommand{\Htwod}[2]{
\begin{tikzpicture}[baseline={([yshift=-.5ex]current bounding box.center)}]
\node[fill=black, circle, inner sep=1pt, minimum size=0.1cm] (1) {};
\node[fill=black, circle, inner sep=1pt, minimum size=0.1cm] (2) [below left = 0.15cm of 1] {};
\node[fill=black, circle, inner sep=1pt, minimum size=0.1cm] (3) [below right =0.15cm of 1] {};
\node[fill=black, circle, inner sep=1pt, minimum size=0.1cm, label=left:\footnotesize{#1}] (4) [below left = 0.15cm of 2] {};
\node[fill=black, circle, inner sep=1pt, minimum size=0.1cm, label=right:\footnotesize{#2}] (5) [below right = 0.15cm of 3] {};
\draw (5)--(3)--(1)--(2)--(4);
\end{tikzpicture}}

\newcommand{\Hthree}[3]{
\begin{tikzpicture}[baseline={([yshift=-.5ex]current bounding box.center)}]
\node[fill=black, circle, inner sep=1pt, minimum size=0.1cm, label=right:\footnotesize{#1}] (1) {};
\node[fill=black, circle, inner sep=1pt, minimum size=0.1cm, label=left:\footnotesize{#2}] (2) [below left = 0.3cm of 1] {};
\node[fill=black, circle, inner sep=1pt, minimum size=0.1cm, label=right:\footnotesize{#3}] (3) [below right = 0.3cm of 1] {};
\draw (1)--(3)--(2)--(1);
\end{tikzpicture}}

\newcommand{\uHthree}[0]{
\begin{tikzpicture}[baseline={([yshift=-.5ex]current bounding box.center)}]
\node[fill=black, circle, inner sep=1pt, minimum size=0.1cm] (1) {};
\node[fill=black, circle, inner sep=1pt, minimum size=0.1cm] (2) [below left = 0.3cm of 1] {};
\node[fill=black, circle, inner sep=1pt, minimum size=0.1cm] (3) [below right = 0.3cm of 1] {};
\draw (1)--(3)--(2)--(1);
\end{tikzpicture}}

\newcommand{\Kfour}[4]{
\begin{tikzpicture}[baseline={([yshift=-.5ex]current bounding box.center)}]
\node[fill=black, circle, inner sep=1pt, minimum size=0.1cm, label=left:\footnotesize{#1}] (1) {};
\node[fill=black, circle, inner sep=1pt, minimum size=0.1cm, label=right:\footnotesize{#2}] (2) [right = 0.3cm of 1] {};
\node[fill=black, circle, inner sep=1pt, minimum size=0.1cm, label=right:\footnotesize{#3}] (3) [below = 0.3cm of 2] {};
\node[fill=black, circle, inner sep=1pt, minimum size=0.1cm, label=left:\footnotesize{#4}] (4) [below = 0.3cm of 1] {};
\draw (4)--(2)--(1)--(3)--(4)--(1)--(2)--(3);
\end{tikzpicture}}

\newcommand{\uKfour}[0]{
\begin{tikzpicture}[baseline={([yshift=-.5ex]current bounding box.center)}]
\node[fill=black, circle, inner sep=1pt, minimum size=0.1cm] (1) {};
\node[fill=black, circle, inner sep=1pt, minimum size=0.1cm] (2) [right = 0.3cm of 1] {};
\node[fill=black, circle, inner sep=1pt, minimum size=0.1cm] (3) [below = 0.3cm of 2] {};
\node[fill=black, circle, inner sep=1pt, minimum size=0.1cm] (4) [below = 0.3cm of 1] {};
\draw (4)--(2)--(1)--(3)--(4)--(1)--(2)--(3);
\end{tikzpicture}}

\newcommand{\threestar}[4]{
\begin{tikzpicture}[baseline={([yshift=-.5ex]current bounding box.center)}]
\node[fill=black, circle, inner sep=1pt, minimum size=0.1cm, label=right:\footnotesize{#1}] (1) {};
\node[fill=black, circle, inner sep=1pt, minimum size=0.1cm, label=left:\footnotesize{#2}] (2) [above = 0.3cm of 1] {};
\node[fill=black, circle, inner sep=1pt, minimum size=0.1cm, label=left:\footnotesize{#3}] (3) [below left = 0.3cm of 1] {};
\node[fill=black, circle, inner sep=1pt, minimum size=0.1cm,label=right:\footnotesize{#4}] (4) [below right = 0.3cm of 1] {};
\draw (2)--(1)--(3)--(1)--(4);
\end{tikzpicture}}

\newcommand{\uthreestar}[0]{
\begin{tikzpicture}[baseline={([yshift=-.5ex]current bounding box.center)}]
\node[fill=black, circle, inner sep=1pt, minimum size=0.1cm] (1) {};
\node[fill=black, circle, inner sep=1pt, minimum size=0.1cm] (2) [above = 0.3cm of 1] {};
\node[fill=black, circle, inner sep=1pt, minimum size=0.1cm] (3) [below left = 0.3cm of 1] {};
\node[fill=black, circle, inner sep=1pt, minimum size=0.1cm] (4) [below right = 0.3cm of 1] {};
\draw (2)--(1)--(3)--(1)--(4);
\end{tikzpicture}}

\newcommand{\usquareteletubby}[0]{
\begin{tikzpicture}[baseline={([yshift=-.5ex]current bounding box.center)}]
\node[fill=black, circle, inner sep=1pt, minimum size=0.1cm] (1) {};
\node[fill=black, circle, inner sep=1pt, minimum size=0.1cm] (2) [above = 0.3cm of 1] {};
\node[fill=black, circle, inner sep=1pt, minimum size=0.1cm] (3) [below left = 0.3cm of 1] {};
\node[fill=black, circle, inner sep=1pt, minimum size=0.1cm] (4) [below right = 0.3cm of 1] {};
\node[fill=black, circle, inner sep=1pt, minimum size=0.1cm] (5) [below left = 0.3cm of 4] {};
\draw (2)--(1)--(3)--(1)--(4)--(5)--(3);
\end{tikzpicture}}

\newcommand{\vedge}[2]{
\begin{tikzpicture}[baseline={([yshift=-.5ex]current bounding box.center)}]
\node[fill=black, circle, inner sep=1pt, minimum size=0.1cm,label=right:\footnotesize{#1}] (1) {};
\node[fill=black, circle, inner sep=1pt, minimum size=0.1cm,label=right:\footnotesize{#2}] (2) [below = 0.3cm of 1] {};
\draw (2)--(1);
\end{tikzpicture}}

\newcommand{\uvedge}[0]{
\begin{tikzpicture}[baseline={([yshift=-.5ex]current bounding box.center)}]
\node[fill=black, circle, inner sep=1pt, minimum size=0.1cm] (1) {};
\node[fill=black, circle, inner sep=1pt, minimum size=0.1cm] (2) [below = 0.3cm of 1] {};
\node[draw=none, fill=none, circle, inner sep=1pt, minimum size=0.1cm] (4) [left = 0.05cm of 1] {};
\draw (2)--(1);
\end{tikzpicture}}

\newcommand{\vnonedge}[2]{
\begin{tikzpicture}[baseline={([yshift=-.5ex]current bounding box.center)}]
\node[fill=black, circle, inner sep=1pt, minimum size=0.1cm, label=right:\footnotesize{#1}] (1) {};
\node[fill=black, circle, inner sep=1pt, minimum size=0.1cm,label=right:\footnotesize{#2}] (2) [below = 0.3cm of 1] {};
\end{tikzpicture}}

\newcommand{\uvedgesquared}[0]{
\begin{tikzpicture}[baseline={([yshift=-.5ex]current bounding box.center)}]
\node[fill=black, circle, inner sep=1pt, minimum size=0.1cm] (1) {};
\node[fill=black, circle, inner sep=1pt, minimum size=0.1cm] (2) [below = 0.3cm of 1] {};
\node[fill=black, circle, inner sep=1pt, minimum size=0.1cm] (3) [right = 0.1cm of 1] {};
\node[fill=black, circle, inner sep=1pt, minimum size=0.1cm] (4) [below = 0.3cm of 3] {};
\draw (2)--(1);
\draw (3)--(4);
\end{tikzpicture}}

\newcommand{\vedgesquared}[4]{
\begin{tikzpicture}[baseline={([yshift=-.5ex]current bounding box.center)}]
\node[fill=black, circle, inner sep=1pt, minimum size=0.1cm,label=left:\footnotesize{#1}] (1) {};
\node[fill=black, circle, inner sep=1pt, minimum size=0.1cm,label=left:\footnotesize{#2}] (2) [below = 0.3cm of 1] {};
\node[fill=black, circle, inner sep=1pt, minimum size=0.1cm,label=right:\footnotesize{#3}] (3) [right = 0.1cm of 1] {};
\node[fill=black, circle, inner sep=1pt, minimum size=0.1cm,label=right:\footnotesize{#4}] (4) [below = 0.3cm of 3] {};
\draw (2)--(1);
\draw (3)--(4);
\end{tikzpicture}}

\newcommand{\Htwosquared}[6]{
\begin{tikzpicture}[baseline={([yshift=-.5ex]current bounding box.center)}]
\node[fill=black, circle, inner sep=1pt, minimum size=0.1cm, label=right:\footnotesize{#1}] (1) {};
\node[fill=black, circle, inner sep=1pt, minimum size=0.1cm, label=left:\footnotesize{#2}] (2) [below left = 0.3cm of 1] {};
\node[fill=black, circle, inner sep=1pt, minimum size=0.1cm, label=right:\footnotesize{#3}] (3) [below right =0.3cm of 1] {};
\node[fill=black, circle, inner sep=1pt, minimum size=0.1cm, label=right:\footnotesize{#4}] (4) [right = 0.65cm of 1] {};
\node[fill=black, circle, inner sep=1pt, minimum size=0.1cm, label=left:\footnotesize{#5}] (5) [below left = 0.3cm of 4] {};
\node[fill=black, circle, inner sep=1pt, minimum size=0.1cm, label=right:\footnotesize{#6}] (6) [below right =0.3cm of 4] {};
\draw (3)--(1)--(2);
\draw (6)--(4)--(5);
\end{tikzpicture}}

\newcommand{\uHtwosquared}[0]{
\begin{tikzpicture}[baseline={([yshift=-.5ex]current bounding box.center)}]
\node[fill=black, circle, inner sep=1pt, minimum size=0.1cm] (1) {};
\node[fill=black, circle, inner sep=1pt, minimum size=0.1cm] (2) [below left = 0.3cm of 1] {};
\node[fill=black, circle, inner sep=1pt, minimum size=0.1cm] (3) [below right =0.3cm of 1] {};
\node[fill=black, circle, inner sep=1pt, minimum size=0.1cm] (4) [right = 0.65cm of 1] {};
\node[fill=black, circle, inner sep=1pt, minimum size=0.1cm] (5) [below left = 0.3cm of 4] {};
\node[fill=black, circle, inner sep=1pt, minimum size=0.1cm] (6) [below right =0.3cm of 4] {};
\draw (3)--(1)--(2);
\draw (6)--(4)--(5);
\end{tikzpicture}}
\newcommand{\fourstar}[5]{
\begin{tikzpicture}[baseline={([yshift=-.5ex]current bounding box.center)}]
\node[fill=black, circle, inner sep=1pt, minimum size=0.1cm, label=left:\footnotesize{#1}] (1) {};
\node[fill=black, circle, inner sep=1pt, minimum size=0.1cm, label=left:\footnotesize{#2}] (2) [above left = 0.3cm of 1] {};
\node[fill=black, circle, inner sep=1pt, minimum size=0.1cm, label=right:\footnotesize{#3}] (3) [above right = 0.3cm of 1] {};
\node[fill=black, circle, inner sep=1pt, minimum size=0.1cm,label=right:\footnotesize{#4}] (4) [below right = 0.3cm of 1] {};
\node[fill=black, circle, inner sep=1pt, minimum size=0.1cm,label=left:\footnotesize{#5}] (5) [below left = 0.3cm of 1] {};
\draw (2)--(1)--(3);
\draw (4)--(1)--(5);
\end{tikzpicture}}

\newcommand{\ufourstar}[0]{
\begin{tikzpicture}[baseline={([yshift=-.5ex]current bounding box.center)}]
\node[fill=black, circle, inner sep=1pt, minimum size=0.1cm] (1) {};
\node[fill=black, circle, inner sep=1pt, minimum size=0.1cm] (2) [above left = 0.3cm of 1] {};
\node[fill=black, circle, inner sep=1pt, minimum size=0.1cm] (3) [above right = 0.3cm of 1] {};
\node[fill=black, circle, inner sep=1pt, minimum size=0.1cm] (4) [below right = 0.3cm of 1] {};
\node[fill=black, circle, inner sep=1pt, minimum size=0.1cm] (5) [below left = 0.3cm of 1] {};
\draw (2)--(1)--(3);
\draw (4)--(1)--(5);
\end{tikzpicture}}

\newcommand{\broom}[5]{
\begin{tikzpicture}[baseline={([yshift=-.5ex]current bounding box.center)}]
\node[fill=black, circle, inner sep=1pt, minimum size=0.1cm, label=left:\footnotesize{#1}] (1) {};
\node[fill=black, circle, inner sep=1pt, minimum size=0.1cm, label=left:\footnotesize{#2}] (2) [below = 0.3cm of 1] {};
\node[fill=black, circle, inner sep=1pt, minimum size=0.1cm, label=left:\footnotesize{#3}] (3) [below = 0.3cm of 2] {};
\node[fill=black, circle, inner sep=1pt, minimum size=0.1cm, label=left:\footnotesize{#5}] (4) [below left = 0.3cm of 3] {};
\node[fill=black, circle, inner sep=1pt, minimum size=0.1cm, label=right:\footnotesize{#4}] (5) [below right = 0.3cm of 3] {};
\draw (1)--(2)--(3);
\draw (4)--(3)--(5);
\end{tikzpicture}}

\newcommand{\ubroom}[0]{
\begin{tikzpicture}[baseline={([yshift=-.5ex]current bounding box.center)}]
\node[fill=black, circle, inner sep=1pt, minimum size=0.1cm] (1) {};
\node[fill=black, circle, inner sep=1pt, minimum size=0.1cm] (2) [below = 0.3cm of 1] {};
\node[fill=black, circle, inner sep=1pt, minimum size=0.1cm] (3) [below = 0.3cm of 2] {};
\node[fill=black, circle, inner sep=1pt, minimum size=0.1cm] (4) [below left = 0.3cm of 3] {};
\node[fill=black, circle, inner sep=1pt, minimum size=0.1cm] (5) [below right = 0.3cm of 3] {};
\draw (1)--(2)--(3);
\draw (4)--(3)--(5);
\end{tikzpicture}}

\newcommand{\upthree}[0]{
\begin{tikzpicture}[baseline={([yshift=-.5ex]current bounding box.center)}]
\node[fill=black, circle, inner sep=1pt, minimum size=0.1cm] (1) {};
\node[fill=black, circle, inner sep=1pt, minimum size=0.1cm] (2) [above = 0.3cm of 1] {};
\node[fill=black, circle, inner sep=1pt, minimum size=0.1cm] (3) [right = 0.3cm of 2] {};
\node[fill=black, circle, inner sep=1pt, minimum size=0.1cm] (4) [below = 0.3cm of 3] {};
\draw (1)--(2)--(3)--(4);
\end{tikzpicture}}

\newcommand{\usquare}[0]{
\begin{tikzpicture}[baseline={([yshift=-.5ex]current bounding box.center)}]
\node[fill=black, circle, inner sep=1pt, minimum size=0.1cm] (1) {};
\node[fill=black, circle, inner sep=1pt, minimum size=0.1cm] (2) [above = 0.3cm of 1] {};
\node[fill=black, circle, inner sep=1pt, minimum size=0.1cm] (3) [right = 0.3cm of 2] {};
\node[fill=black, circle, inner sep=1pt, minimum size=0.1cm] (4) [below = 0.3cm of 3] {};
\draw (1)--(2)--(3)--(4)--(1);
\end{tikzpicture}}

\newcommand{\pthree}[4]{
\begin{tikzpicture}[baseline={([yshift=-.5ex]current bounding box.center)}]
\node[fill=black, circle, inner sep=1pt, minimum size=0.1cm, label=left:\footnotesize{#1}] (1) {};
\node[fill=black, circle, inner sep=1pt, minimum size=0.1cm, label=left:\footnotesize{#2}] (2) [above = 0.3cm of 1] {};
\node[fill=black, circle, inner sep=1pt, minimum size=0.1cm, label=right:\footnotesize{#3}] (3) [right = 0.3cm of 2] {};
\node[fill=black, circle, inner sep=1pt, minimum size=0.1cm, label=right:\footnotesize{#4}] (4) [below = 0.3cm of 3] {};
\draw (1)--(2)--(3)--(4);
\end{tikzpicture}}

\newcommand{\ucfour}[0]{
\begin{tikzpicture}[baseline={([yshift=-.5ex]current bounding box.center)}]
\node[fill=black, circle, inner sep=1pt, minimum size=0.1cm] (1) {};
\node[fill=black, circle, inner sep=1pt, minimum size=0.1cm] (2) [above = 0.3cm of 1] {};
\node[fill=black, circle, inner sep=1pt, minimum size=0.1cm] (3) [right = 0.3cm of 2] {};
\node[fill=black, circle, inner sep=1pt, minimum size=0.1cm] (4) [below = 0.3cm of 3] {};
\draw (1)--(2)--(3)--(4)--(1);
\end{tikzpicture}}

\newcommand{\uceight}[0]{
\begin{tikzpicture}[baseline={([yshift=-.5ex]current bounding box.center)}]
\node[fill=black, circle, inner sep=1pt, minimum size=0.1cm] (1) {};
\node[fill=black, circle, inner sep=1pt, minimum size=0.1cm] (2) [above = 0.15cm of 1] {};
\node[fill=black, circle, inner sep=1pt, minimum size=0.1cm] (3) [above = 0.15cm of 2] {};
\node[fill=black, circle, inner sep=1pt, minimum size=0.1cm] (4) [right = 0.15cm of 3] {};
\node[fill=black, circle, inner sep=1pt, minimum size=0.1cm] (5) [right = 0.15cm of 4] {};
\node[fill=black, circle, inner sep=1pt, minimum size=0.1cm] (6) [below = 0.15cm of 5] {};
\node[fill=black, circle, inner sep=1pt, minimum size=0.1cm] (7) [below = 0.15cm of 6] {};
\node[fill=black, circle, inner sep=1pt, minimum size=0.1cm] (8) [left = 0.15cm of 7] {};
\draw (1)--(2)--(3)--(4)--(5)--(6)--(7)--(8)--(1);
\end{tikzpicture}}

\newcommand{\upfour}[0]{
\begin{tikzpicture}[baseline={([yshift=-.5ex]current bounding box.center)}]
\node[fill=black, circle, inner sep=1pt, minimum size=0.1cm] (1) {};
\node[fill=black, circle, inner sep=1pt, minimum size=0.1cm] (2) [above = 0.3cm of 1] {};
\node[fill=black, circle, inner sep=1pt, minimum size=0.1cm] (3) [right = 0.3cm of 2] {};
\node[fill=black, circle, inner sep=1pt, minimum size=0.1cm] (4) [below = 0.3cm of 3] {};
\node[fill=black, circle, inner sep=1pt, minimum size=0.1cm] (5) [right = 0.3cm of 4] {};
\draw (1)--(2)--(3)--(4)--(5);
\end{tikzpicture}}

\newcommand{\pfour}[5]{
\begin{tikzpicture}[baseline={([yshift=-.5ex]current bounding box.center)}]
\node[fill=black, circle, inner sep=1pt, minimum size=0.1cm, label=below:\footnotesize{#1}] (1) {};
\node[fill=black, circle, inner sep=1pt, minimum size=0.1cm, label=above:\footnotesize{#2}] (2) [above = 0.3cm of 1] {};
\node[fill=black, circle, inner sep=1pt, minimum size=0.1cm, label=above:\footnotesize{#3}] (3) [right = 0.3cm of 2] {};
\node[fill=black, circle, inner sep=1pt, minimum size=0.1cm, label=below:\footnotesize{#4}] (4) [below = 0.3cm of 3] {};
\node[fill=black, circle, inner sep=1pt, minimum size=0.1cm, label=below:\footnotesize{#5}] (5) [right = 0.3cm of 4] {};
\draw (1)--(2)--(3)--(4)--(5);
\end{tikzpicture}}

\newcommand{\upfive}[0]{
\begin{tikzpicture}[baseline={([yshift=-.5ex]current bounding box.center)}]
\node[fill=black, circle, inner sep=1pt, minimum size=0.1cm] (1) {};
\node[fill=black, circle, inner sep=1pt, minimum size=0.1cm] (2) [above = 0.3cm of 1] {};
\node[fill=black, circle, inner sep=1pt, minimum size=0.1cm] (3) [right = 0.3cm of 2] {};
\node[fill=black, circle, inner sep=1pt, minimum size=0.1cm] (4) [below = 0.3cm of 3] {};
\node[fill=black, circle, inner sep=1pt, minimum size=0.1cm] (5) [right = 0.3cm of 4] {};
\node[fill=black, circle, inner sep=1pt, minimum size=0.1cm] (6) [above = 0.3cm of 5] {};
\draw (1)--(2)--(3)--(4)--(5)--(6);
\end{tikzpicture}}

\newcommand{\upseven}[0]{
\begin{tikzpicture}[baseline={([yshift=-.5ex]current bounding box.center)}]
\node[fill=black, circle, inner sep=1pt, minimum size=0.1cm] (1) {};
\node[fill=black, circle, inner sep=1pt, minimum size=0.1cm] (2) [above = 0.3cm of 1] {};
\node[fill=black, circle, inner sep=1pt, minimum size=0.1cm] (3) [right = 0.3cm of 2] {};
\node[fill=black, circle, inner sep=1pt, minimum size=0.1cm] (4) [below = 0.3cm of 3] {};
\node[fill=black, circle, inner sep=1pt, minimum size=0.1cm] (5) [right = 0.3cm of 4] {};
\node[fill=black, circle, inner sep=1pt, minimum size=0.1cm] (6) [above = 0.3cm of 5] {};
\node[fill=black, circle, inner sep=1pt, minimum size=0.1cm] (7) [right = 0.3cm of 6] {};
\node[fill=black, circle, inner sep=1pt, minimum size=0.1cm] (8) [below = 0.3cm of 7] {};
\draw (1)--(2)--(3)--(4)--(5)--(6)--(7)--(8);
\end{tikzpicture}}

\newcommand{\uyseven}[0]{
\begin{tikzpicture}[baseline={([yshift=-.5ex]current bounding box.center)}]
\node[fill=black, circle, inner sep=1pt, minimum size=0.1cm] (1) {};
\node[fill=black, circle, inner sep=1pt, minimum size=0.1cm] (2) [above = 0.3cm of 1] {};
\node[fill=black, circle, inner sep=1pt, minimum size=0.1cm] (3) [right = 0.3cm of 2] {};
\node[fill=black, circle, inner sep=1pt, minimum size=0.1cm] (4) [below = 0.3cm of 3] {};
\node[fill=black, circle, inner sep=1pt, minimum size=0.1cm] (5) [right = 0.3cm of 4] {};
\node[fill=black, circle, inner sep=1pt, minimum size=0.1cm] (6) [above = 0.3cm of 5] {};
\node[fill=black, circle, inner sep=1pt, minimum size=0.1cm] (7) [right = 0.3cm of 6] {};
\node[fill=black, circle, inner sep=1pt, minimum size=0.1cm] (8) [above = 0.3cm of 3] {};
\draw (1)--(2)--(3)--(4)--(5)--(6)--(7);
\draw (3)--(8);
\end{tikzpicture}}

\newcommand{\usquareteletubbytwoantennas}[0]{
\begin{tikzpicture}[baseline={([yshift=-.5ex]current bounding box.center)}]
\node[fill=black, circle, inner sep=1pt, minimum size=0.1cm] (1) {};
\node[fill=black, circle, inner sep=1pt, minimum size=0.1cm] (2) [above = 0.3cm of 1] {};
\node[fill=black, circle, inner sep=1pt, minimum size=0.1cm] (3) [right = 0.3cm of 2] {};
\node[fill=black, circle, inner sep=1pt, minimum size=0.1cm] (4) [below = 0.3cm of 3] {};
\node[fill=black, circle, inner sep=1pt, minimum size=0.1cm] (5) [right = 0.3cm of 4] {};
\node[fill=black, circle, inner sep=1pt, minimum size=0.1cm] (6) [above = 0.3cm of 5] {};
\draw (1)--(2)--(3)--(4)--(1);
\draw (4)--(5);
\draw (3)--(6);
\end{tikzpicture}}

\newcommand{\ustaronetwothree}[0]{
\begin{tikzpicture}[baseline={([yshift=-.5ex]current bounding box.center)}]
\node[fill=black, circle, inner sep=1pt, minimum size=0.1cm] (1) {};
\node[fill=black, circle, inner sep=1pt, minimum size=0.1cm] (2) [above = 0.3cm of 1] {};
\node[fill=black, circle, inner sep=1pt, minimum size=0.1cm] (3) [right = 0.3cm of 2] {};
\node[fill=black, circle, inner sep=1pt, minimum size=0.1cm] (4) [below = 0.3cm of 3] {};
\node[fill=black, circle, inner sep=1pt, minimum size=0.1cm] (5) [right = 0.3cm of 4] {};
\node[fill=black, circle, inner sep=1pt, minimum size=0.1cm] (6) [above = 0.3cm of 5] {};
\node[fill=black, circle, inner sep=1pt, minimum size=0.1cm] (7) [right = 0.3cm of 5] {};
\draw (1)--(4)--(5)--(7);
\draw (2)--(3)--(6)--(5);
\end{tikzpicture}}

\newcommand{\uuneventeletubby}[0]{
\begin{tikzpicture}[baseline={([yshift=-.5ex]current bounding box.center)}]
\node[fill=black, circle, inner sep=1pt, minimum size=0.1cm] (1) {};
\node[fill=black, circle, inner sep=1pt, minimum size=0.1cm] (2) [above = 0.3cm of 1] {};
\node[fill=black, circle, inner sep=1pt, minimum size=0.1cm] (3) [right = 0.3cm of 2] {};
\node[fill=black, circle, inner sep=1pt, minimum size=0.1cm] (4) [below = 0.3cm of 3] {};
\node[fill=black, circle, inner sep=1pt, minimum size=0.1cm] (5) [right = 0.3cm of 4] {};
\node[fill=black, circle, inner sep=1pt, minimum size=0.1cm] (6) [above = 0.3cm of 5] {};
\node[fill=black, circle, inner sep=1pt, minimum size=0.1cm] (7) [right = 0.3cm of 6] {};
\draw (1)--(2)--(3)--(4)--(1);
\draw (4)--(5);
\draw (3)--(6)--(7);
\end{tikzpicture}}

\newcommand{\ufivestar}[0]{
\begin{tikzpicture}[baseline={([yshift=-.5ex]current bounding box.center)}]
\node[fill=black, circle, inner sep=1pt, minimum size=0.1cm] (1) at (0,0) {};
\node[fill=black, circle, inner sep=1pt, minimum size=0.1cm] (2) at ($(1)+(0:0.3cm)$) {};
\node[fill=black, circle, inner sep=1pt, minimum size=0.1cm] (3) at ($(1)+(72:0.3cm)$) {};
\node[fill=black, circle, inner sep=1pt, minimum size=0.1cm] (4) at ($(1)+(144:0.3cm)$) {};
\node[fill=black, circle, inner sep=1pt, minimum size=0.1cm] (5) at ($(1)+(216:0.3cm)$) {};
\node[fill=black, circle, inner sep=1pt, minimum size=0.1cm] (6) at ($(1)+(288:0.3cm)$) {};
\draw (2)--(1)--(3);
\draw (4)--(1)--(5);
\draw (1)--(6);
\end{tikzpicture}}

\newcommand{\uweights}[0]{
\begin{tikzpicture}[baseline={([yshift=-.5ex]current bounding box.center)}]
\node[fill=black, circle, inner sep=1pt, minimum size=0.1cm, label] (1) {};
\node[fill=black, circle, inner sep=1pt, minimum size=0.1cm, label] (2) [above left= 0.3cm of 1] {};
\node[fill=black, circle, inner sep=1pt, minimum size=0.1cm, label] (3) [below left = 0.3cm of 1] {};
\node[fill=black, circle, inner sep=1pt, minimum size=0.1cm,label] (4) [right = 0.3cm of 1] {};
\node[fill=black, circle, inner sep=1pt, minimum size=0.1cm,label] (5) [above right = 0.3cm of 4] {};
\node[fill=black, circle, inner sep=1pt, minimum size=0.1cm,label] (6) [below right = 0.3cm of 4] {};
\draw (2)--(1)--(3)--(1)--(4)--(5);
\draw (4)--(6);
\end{tikzpicture}}

\newcommand{\uspider}[0]{
\begin{tikzpicture}[baseline={([yshift=-.5ex]current bounding box.center)}]
\node[fill=black, circle, inner sep=1pt, minimum size=0.1cm, label] (1) {};
\node[fill=black, circle, inner sep=1pt, minimum size=0.1cm, label] (2) [above left= 0.3cm of 1] {};
\node[fill=black, circle, inner sep=1pt, minimum size=0.1cm, label] (3) [below left = 0.3cm of 1] {};
\node[fill=black, circle, inner sep=1pt, minimum size=0.1cm, label] (7) [left = 0.3cm of 3] {};
\node[fill=black, circle, inner sep=1pt, minimum size=0.1cm,label] (4) [right = 0.3cm of 1] {};
\node[fill=black, circle, inner sep=1pt, minimum size=0.1cm,label] (5) [above right = 0.3cm of 4] {};
\node[fill=black, circle, inner sep=1pt, minimum size=0.1cm,label] (6) [below right = 0.3cm of 4] {};
\node[fill=black, circle, inner sep=1pt, minimum size=0.1cm,label] (8) [right = 0.3cm of 6] {};
\draw (2)--(1)--(3)--(1)--(4)--(5);
\draw (4)--(6)--(8);
\draw (3)--(7);
\end{tikzpicture}}

\newcommand{\uloongweights}[0]{
\begin{tikzpicture}[baseline={([yshift=-.5ex]current bounding box.center)}]
\node[fill=black, circle, inner sep=1pt, minimum size=0.1cm, label] (1) {};
\node[fill=black, circle, inner sep=1pt, minimum size=0.1cm, label] (2) [above left= 0.3cm of 1] {};
\node[fill=black, circle, inner sep=1pt, minimum size=0.1cm, label] (3) [below left = 0.3cm of 1] {};
\node[fill=black, circle, inner sep=1pt, minimum size=0.1cm,label] (4) [right = 0.3cm of 1] {};
\node[fill=black, circle, inner sep=1pt, minimum size=0.1cm,label] (5) [right = 0.3cm of 4] {};
\node[fill=black, circle, inner sep=1pt, minimum size=0.1cm,label] (6) [right = 0.3cm of 5] {};
\node[fill=black, circle, inner sep=1pt, minimum size=0.1cm,label] (7) [above right = 0.3cm of 6] {};
\node[fill=black, circle, inner sep=1pt, minimum size=0.1cm,label] (8) [below right = 0.3cm of 6] {};
\draw (2)--(1)--(3)--(1)--(4)--(5)--(6)--(7);
\draw (6)--(8);
\end{tikzpicture}}

\newcommand{\ulongbroom}[0]{
\begin{tikzpicture}[baseline={([yshift=-.5ex]current bounding box.center)}]
\node[fill=black, circle, inner sep=1pt, minimum size=0.1cm, label] (1) {};
\node[fill=black, circle, inner sep=1pt, minimum size=0.1cm, label] (2) [above left= 0.3cm of 1] {};
\node[fill=black, circle, inner sep=1pt, minimum size=0.1cm, label] (3) [below left = 0.3cm of 1] {};
\node[fill=black, circle, inner sep=1pt, minimum size=0.1cm,label] (4) [right = 0.3cm of 1] {};
\node[fill=black, circle, inner sep=1pt, minimum size=0.1cm,label] (5) [right = 0.3cm of 4] {};
\node[fill=black, circle, inner sep=1pt, minimum size=0.1cm,label] (6) [right = 0.3cm of 5] {};
\draw (2)--(1)--(3)--(1)--(4)--(5)--(6);
\end{tikzpicture}}

\newcommand{\uthreebroom}[0]{
\begin{tikzpicture}[baseline={([yshift=-.5ex]current bounding box.center)}]
\node[fill=black, circle, inner sep=1pt, minimum size=0.1cm, label] (1) {};
\node[fill=black, circle, inner sep=1pt, minimum size=0.1cm, label] (2) [above left= 0.3cm of 1] {};
\node[fill=black, circle, inner sep=1pt, minimum size=0.1cm, label] (3) [below left = 0.3cm of 1] {};
\node[fill=black, circle, inner sep=1pt, minimum size=0.1cm,label] (4) [right = 0.3cm of 1] {};
\node[fill=black, circle, inner sep=1pt, minimum size=0.1cm,label] (5) [right = 0.3cm of 4] {};
\node[fill=black, circle, inner sep=1pt, minimum size=0.1cm,label] (6) [left = 0.3cm of 1] {};
\draw (2)--(1)--(3)--(1)--(4)--(5);
\draw (1)--(6);
\end{tikzpicture}}

\newcommand{\uchristmastree}[0]{
\begin{tikzpicture}[baseline={([yshift=-.5ex]current bounding box.center)}]
\node[fill=black, circle, inner sep=1pt, minimum size=0.1cm, label] (1) {};
\node[fill=black, circle, inner sep=1pt, minimum size=0.1cm, label] (2) [left= 0.3cm of 1] {};
\node[fill=black, circle, inner sep=1pt, minimum size=0.1cm, label] (3) [below right= 0.3cm of 1] {};
\node[fill=black, circle, inner sep=1pt, minimum size=0.1cm,label] (4) [above right = 0.3cm of 1] {};
\node[fill=black, circle, inner sep=1pt, minimum size=0.1cm,label] (5) [right = 0.3cm of 4] {};
\node[fill=black, circle, inner sep=1pt, minimum size=0.1cm,label] (6) [right = 0.3cm of 3] {};
\draw (2)--(1)--(3)--(1)--(4)--(5);
\draw (3)--(6);
\end{tikzpicture}}

\newcommand{\threestarplus}[4]{
\begin{tikzpicture}[baseline={([yshift=-.5ex]current bounding box.center)}]
\node[fill=black, circle, inner sep=1pt, minimum size=0.1cm, label=right:\footnotesize{#1}] (1) {};
\node[fill=black, circle, inner sep=1pt, minimum size=0.1cm, label=left:\footnotesize{#2}] (2) [above = 0.3cm of 1] {};
\node[fill=black, circle, inner sep=1pt, minimum size=0.1cm, label=left:\footnotesize{#3}] (3) [below left = 0.3cm of 1] {};
\node[fill=black, circle, inner sep=1pt, minimum size=0.1cm,label=right:\footnotesize{#4}] (4) [below right = 0.3cm of 1] {};
\draw (2)--(1)--(3)--(1)--(4);
\draw (3)--(4);
\end{tikzpicture}}

\newcommand{\bottomlesshouse}[5]{
\begin{tikzpicture}[baseline={([yshift=-.5ex]current bounding box.center)}]
\node[fill=black, circle, inner sep=1pt, minimum size=0.1cm, label=right:\footnotesize{#1}] (1) {};
\node[fill=black, circle, inner sep=1pt, minimum size=0.1cm, label=left:\footnotesize{#2}] (2) [below left = 0.3cm of 1] {};
\node[fill=black, circle, inner sep=1pt, minimum size=0.1cm, label=right:\footnotesize{#3}] (3) [below right = 0.3cm of 1] {};
\node[fill=black, circle, inner sep=1pt, minimum size=0.1cm, label=left:\footnotesize{#4}] (4) [below = 0.3cm of 2] {};
\node[fill=black, circle, inner sep=1pt, minimum size=0.1cm, label=right:\footnotesize{#5}] (5) [below = 0.3cm of 3] {};
\draw (1)--(3)--(2)--(1);
\draw (2)--(4);
\draw (3)--(5);
\end{tikzpicture}}

\newcommand{\Kfivedconstruction}[7]{
\begin{tikzpicture}[baseline={([yshift=-.5ex]current bounding box.center)}]
\node[fill=black, circle, inner sep=1pt, minimum size=0.1cm, label={[label distance=0.001cm]below:\footnotesize{#1}}] (1) {};
\node[fill=black, circle, inner sep=1pt, minimum size=0.1cm] (2) [above = 0.5cm of 1] {};
\node[fill=black, circle, inner sep=1pt, minimum size=0.1cm] (3) [below left = 0.71cm of 2] {};
\node[fill=black, circle, inner sep=1pt, minimum size=0.1cm] (4) [below right = 0.71cm of 2] {};
\node[fill=black, circle, inner sep=1pt, minimum size=0.1cm, label=below:\footnotesize{#2}] (5) [below left = 0.71cm of 3] {};
\node[fill=black, circle, inner sep=1pt, minimum size=0.1cm, label=below:\footnotesize{#3}] (6) [right = 1.05cm of 5] {};
\node[fill=black, circle, inner sep=1pt, minimum size=0.1cm, label=below:\footnotesize{#4}] (7) [below right = 0.71cm of 4] {};
\node[fill=black, circle, inner sep=1pt, minimum size=0.1cm, label=left:\footnotesize{#5}] (8) [left = 0.55cm of 2] {};
\node[fill=black, circle, inner sep=1pt, minimum size=0.1cm, label=right:\footnotesize{#6}] (9) [right = 0.55cm of 2] {};
\node[fill=black, circle, inner sep=1pt, minimum size=0.1cm, label=right:\footnotesize{#7}] (10) [above = 1cm of 2] {};
\node[fill=black, circle, inner sep=1pt, minimum size=0.1cm] (11) [below = 0.5cm of 10] {};
\draw (1)--(2)--(3)--(5)--(6)--(7)--(2);
\draw (5)--(8);
\draw (7)--(9);
\draw (8)--(10)--(9);
\draw (10)--(11)--(2);
\end{tikzpicture}}

\newcommand{\Kfourdconstruction}[4]{
\begin{tikzpicture}[baseline={([yshift=-.5ex]current bounding box.center)}]
\node[fill=black, circle, inner sep=1pt, minimum size=0.1cm, label=right:\footnotesize{#1}] (1) {};
\node[fill=black, circle, inner sep=1pt, minimum size=0.1cm] (2) [below = 0.3 cm of 1] {};
\node[fill=black, circle, inner sep=1pt, minimum size=0.1cm] (3) [below left = 0.3cm of 2] {};
\node[fill=black, circle, inner sep=1pt, minimum size=0.1cm] (4) [below right = 0.3cm of 2] {};
\node[fill=black, circle, inner sep=1pt, minimum size=0.1cm, label=below:\footnotesize{#2}] (5) [below left = 0.3cm of 3] {};
\node[fill=black, circle, inner sep=1pt, minimum size=0.1cm, label=below:\footnotesize{#3}] (6) [right = 0.47cm of 5] {};
\node[fill=black, circle, inner sep=1pt, minimum size=0.1cm, label=below:\footnotesize{#4}] (7) [below right = 0.3cm of 4] {};
\draw (1)--(2)--(3)--(5)--(6)--(7)--(2);
\end{tikzpicture}}

\def \trop {\operatorname{trop}}
\def \hom {\operatorname{hom}}
\def \Hom {\operatorname{Hom}}
\def \HDE {\operatorname{HDE}}
\newcommand{\GU}{{\mathcal{G}_\mathcal{U}}}
\newcommand{\NU}{{\mathcal{N}_\mathcal{U}}}
\newcommand{\RR}{\mathbb{R}}
\newcommand{\ZZ}{\mathbb{Z}}
\newcommand \U {\mathcal{U}}
 \newcommand{\A}{{\mathcal{A}}}

\date{\today}

\title{Non-trivial squares and Sidorenko's conjecture}

\author{Pranav Garg}
\address{Department of Mathematics and Statistics,
Lederle Graduate Research Tower, 1623D,
University of Massachusetts Amherst,
710 N. Pleasant Street,
Amherst, MA 01003} \email{pranavgarg@umass.edu}

\author{Annie Raymond}
\address{Department of Mathematics and Statistics,
Lederle Graduate Research Tower, 1623D,
University of Massachusetts Amherst,
710 N. Pleasant Street,
Amherst, MA 01003} \email{raymond@math.umass.edu}

\author{Amanda Redlich}
\address{Department of Mathematical Sciences,
Southwick Hall 303,
11 University Ave.,
University of Massachusetts Lowell,
Lowell, MA 01854}\email{Amanda\_Redlich@uml.edu}

\thanks{Annie Raymond is partially supported by NSF grant DMS-2054404.}

\begin{document}

\maketitle

\begin{abstract}
    Let $t(H;G)$ be the homomorphism density of a graph $H$ into a graph $G$. Sidorenko's conjecture states that for any bipartite graph $H$, $t(H;G)\geq t(K_2;G)^{|E(H)|}$ for all graphs $G$. It is already known that such inequalities cannot be certified through the sums of squares method when $H$ is a so-called \emph{trivial square}. In this paper, we investigate recent results about Sidorenko's conjecture and classify those involving trivial versus non-trivial squares. We then present some computational results. In particular, we categorize the bipartite graphs $H$ on at most 7 edges for which $t(H;G)\geq t(K_2;G)^{|E(H)|}$ has a sum of squares certificate. We then discuss other limitations for sums of squares proofs beyond trivial squares.
\end{abstract}

\section{Introduction}

One important conjecture in extremal graph theory is Sidorenko's conjecture~(\cite{Sid93}) which has garnered a lot of attention throughout the years (\cite{li2011logarithimic}, \cite{sze}, \cite{conlon2018sidorenko}, \cite{conlonkimleelee}, \cite{KLL}). Many wondered whether the Cauchy-Schwarz method (which is equivalent to the sum of squares method (\cite{blekherman2012semidefinite})) was the right tool to attack this problem (see \cite{gowersblog} for some context). In \cite{BRST}, the authors showed that there are instances of Sidorenko's conjecture that can not be verified with sums of squares. In particular, they showed that it cannot be verified in the case of \emph{trivial squares} (defined below). 

In this paper, we analyze the possibility of sum of squares proofs of Sidorenko's conjecture for multiple graphs and graph families.  We first show that many families of graphs for which Sidorenko's conjecture is known to hold involve non-trivial squares.  Therefore, these could potentially be recovered using sums of squares. 

However, we also show that trivial squares are not the only obstacle for sums of squares certificates. We compute a poset of all bipartite graphs with $k$ edges for $k\leq 7$ ordered by the relation $\succeq$ where $H_1 \succeq H_2$ if $t(H_1;G)-t(H_2;G) \geq 0$ is valid for all graphs $G$, and can be certified with a sum of squares. If Sidorenko's conjecture could be proven through sums of squares, $\uvedge^k$ would be the unique minimal element of that poset. If trivial squares were the only ``important'' obstacle to having sums of squares certificates for Sidorenko's conjecture, all minimal elements would involve trivial squares. However, we demonstrate that there are minimal elements that are not trivial squares, for example, \uloongweights.

\subsection{Preliminaries} 
A graph $G$ has vertex set $V(G)$ and edge set $E(G)$. All graphs are assumed to be simple, without loops or multiple edges. The \emph{homomorphism density} of a graph $H$ in a graph $G$, denoted by $t(H;G)$, is the probability that a random map from $V(H)$ to $V(G)$ is a graph homomorphism, i.e., that it maps every edge of $H$ to an edge of $G$.

We now introduce the gluing algebra of graphs (see \cite{LovaszBook} for a broader exposition). A graph is {\em partially labeled} if a subset $L$ of its vertices are labeled with elements of $\mathbb{N} := \{1,2,3,\ldots\}$ such that no vertex receives more than one label. If no
vertices of $H$ are labeled, then $H$ is {\em unlabeled}. Let $\vartheta: L \rightarrow V(G)$. Then for a partially labeled graph $H$, $t(H;G)$ is the probability that a random map that extends $\vartheta$ is a graph homomorphism.

Let $\A$ denote the vector space of all formal finite $\RR$-linear combinations of partially labeled graphs without isolated vertices, including
the empty graph with no vertices which we denote as $1$. We call an element $a = \sum \alpha_i H_i$ of $\A$ a {\em graph combination}, and each $\alpha_i H_i$ a {\em term} of $a$. 

The vector space $\A$ has a product defined as follows. For two labeled graphs   $H_1$ and $H_2$, form the new labeled graph $H_1H_2$  by gluing together the vertices in the two graphs with the same label, and keeping only one copy of any edge that may have doubled in the process. For example, $\threestarplus{3}{}{1}{2} \cdot \Htwo{1}{2}{}= \bottomlesshouse{2}{3}{1}{}{}$. Equipped with this product, $\A$ becomes an $\RR$-algebra with the empty graph as its multiplicative identity. Note also that any polynomial can be viewed as a graph combination.

For a fixed finite set of labels $L \subset \mathbb{N}$, let $\A_L$ denote the subalgebra of $\A$ spanned by all
graphs whose label sets are contained in $L$. Then $\A_\emptyset$ is the subalgebra of $\A$ spanned by unlabeled graphs. The algebra $\A$ admits a simple linear map into $\A_\emptyset$ that removes the labels in a graph combination to create a graph combination of unlabeled graphs. We call this map {\em unlabeling} and denote it by $[[ \cdot ]]$.  We view elements $a \in \A_\emptyset$ as functions that can be evaluated on unlabeled graphs $G$ via homomorphism densities. An element $a = \sum \alpha_i H_i$ of $\A_\emptyset $ is called nonnegative if $ \sum \alpha_i t(H_i;G)\geq 0$ for all graphs $G$. In this setting, a {\em sum of squares (sos)} in $\A_\emptyset$ is a finite sum of unlabeled squares of graph combinations $a_i \in \A$, namely, $\sum [[ a_i^2 ]]$. Note that unlabeling corresponds to symmetrization, and so this is truly a sum of squares.   It can be easily seen that an sos is a nonnegative graph combination. 

In the language of the gluing algebra, Sidorenko's conjecture states that for every bipartite graph $H$, the graph combination $H-\uvedge^{|E(H)|}$ is nonnegative. Further, for a particular graph $H$, if $H-\uvedge^{|E(H)|}\geq 0$, we say that $H$ is \emph{Sidorenko}.

\begin{example}\label{ex:mult}
One can prove Sidorenko's conjecture when $H=\uHtwo$ with the following sos proof: 

$$ \left[\left[\left(\vedge{1}{\textcolor{white}{1}}-\uvedge\right)^2\right]\right] =  \left[\left[\Htwo{1}{}{}-2\ \vedge{}{} \vedge{1}{\textcolor{white}{1}}+\uvedge\uvedge \ \right]\right]= \uHtwo-\uvedge\ \uvedge \ .$$
\end{example}

It turns out that $H-\uvedge^{|E(H)|}$ is not sos for every bipartite graph $H$. One definition that is important to partially characterize what cases of Sidorenko's conjecture cannot be written as an sos is the following. An unlabeled graph is called a \emph{trivial square} if whenever $H=[[F^2]]$ for some partially labeled graph $F$, then $F$ is a fully labeled copy of $H$. For example, \upthree \  is a trivial square but \uHtwo \  is not since $\uHtwo=\left[\left[\left(\vedge{1}{}\right)^2\right]\right]$.  If a graph $H$ is \emph{not} a trivial square, i.e., there is some $F$ that is not a fully labeled copy of $H$ such that $H=[[F^2]]$, we call $F$ a \emph{square root} of $H$.  Note that $F$ is not necessarily unique; for instance, in Section \ref{sec:previousresultsandsos}, we construct multiple square roots for several classes of graphs.

From \cite{BRST}, we know that if $H$ is not a trivial square, then there must be an involution $\varphi$ of $H$ that fixes a non-empty proper subset of $V(H)$ such that the vertices $v$ that are not fixed by $\varphi$ can be partitioned into two groups, each group consisting of one vertex from each pair $\{v, \varphi(v)\}$, such that there are no edges between vertices in the two groups.  This is Lemma 2.8 in \cite{BRST}.  The proof relies on the effect squaring has on labeled and unlabeled vertices; see the reference for details.

From \cite{BRST}, we also know that if $H$ is a trivial square, then $H-\uvedge^{|E(H)|}$ is not an sos. This is the case for example when $H$ is a path of odd length, the so-called Blakley-Roy inequalities. In that case, Sidorenko's conjecture cannot be proven by the sos method (even though it has been verified with other techniques). This is also the case for the smallest unresolved instance of Sidorenko's conjecture: $K_{5,5}\backslash C_{10}$.  In Sections \ref{sec:computationalresults} and \ref{sec:nosos}, we give additional scenarios in which $H-\uvedge^{|E(H)|}$ is not sos \emph{even though} $H$ itself is a non-trivial square.

\subsection{Organization of the paper} In this paper, we show in Section \ref{sec:previousresultsandsos} that many of the cases for which Sidorenko's conjecture is known to hold involve bipartite graphs that are not trivial squares; we also give constructions of square roots for these graphs. In particular, we show that the recent Cartesian product and subdivision constructions of \cite{conlonkimleelee} and the degree construction of \cite{conlon2018sidorenko} are all non-trivial squares. In the case of the Cartesian product construction in \cite{conlonkimleelee}, we also give a detailed classification of all possible square roots.  The classification and construction of square roots in this way is a key first step in analyzing the potential for an sos proof.

However, non-trivial squares are just one of many obstacles to sos proofs: even if $H$ is a non-trivial square, there is no guarantee that $H-\uvedge^{|E(H)|}$ has a sos proof. In Section \ref{sec:computationalresults}, we look at which graphs $H$ with at most seven edges are such that there is an sos proof for $H-\uvedge^{|E(H)|}$. We then discuss in Section \ref{sec:nosos} some additional road blocks to sums of squares. 


\section{Previous results and sums of squares}\label{sec:previousresultsandsos}

We first go over a few lemmas about non-trivial squares before considering different past results and reframing them within the lens of non-trivial squares.

As explained in the introduction, Lemma 2.8 of \cite{BRST} shows that a \emph{non-trivial square} is a graph $H$ for which there exists an involution $\varphi$ of $H$ that fixes a non-empty proper subset of $V(H)$ such that the vertices $v$ that are not fixed by $\varphi$ can be partitioned into two groups, each group consisting of one vertex from each pair $\{v, \varphi(v)\}$, such that there are no edges between vertices in the two groups. This definition is not always the easiest to work with. We thus note a few simple types of graphs that are always squares.

\begin{lemma}\label{twins}
If there exists two (non-adjacent) vertices $v, v' \in V(H)$ such that $v, v'$ have the same neighborhood, then $H$ is a non-trivial square.
\end{lemma}

\begin{proof}
Let $L$ be the graph $H\backslash{v'}$ where all the vertices are labeled except $v$. Since $[[L^2]]=H$, $H$ is a non-trivial square and $L$ is (one of) its square root(s).  

\end{proof}

In this case, the involution described earlier is simply the $\varphi$ which fixes $H\backslash \{v, v'\}$ and partitions the remainder into $\{v\}$ and $\{v'\}$.  Note that, in fact, this proof approach can be expanded to generate multiple square roots if the graph satisfies certain conditions.  For any $S_1, S_2, \ldots S_k \subseteq V(H)$ such that
\begin{itemize}
    \item the $S_i$ are disjoint sets of even size,
    \item for every $i$, each vertex in $S_i$ has the same neighborhood in $H$, and
    \item no vertex in $\cup_{i=1}^{k}S_i$ neighbors any other vertex in $\cup_{i=1}^{k}S_i$.
\end{itemize}
and for any subsets $S_i' \subset S_i$ such that $2|S_i'|=|S_i|$, the graph $L=H \backslash \left( \cup_{i=1}^{k} S_i'\right)$ with all vertices labeled except those in $\cup_{i=1}^{k} S_i$ is a square root of $H$.  Now the involution is the one which fixes $H \backslash \left( \cup_{i=1}^{k} S_i \right)$ and partitions the remainder into $\cup_{i=1}^{k} S_i'$ and $\cup_{i=1}^{k} S_i \backslash \cup_{i=1}^{k} S_i'$

\begin{lemma}\label{twincc}
If $H$ contains two or more isomorphic connected components, then $H$ is a non-trivial square.
\end{lemma}

\begin{proof}
Let $C_1$ and $C_2$ be isomorphic connected components of $H$ and let $L$ be the graph $H \backslash C_2$ where all the vertices are labeled except those in $C_1$.  Since $[[L^2]]=H$, $H$ is a non-trivial square and $L$ is (one of) its square root(s).
\end{proof}

As before, this construction also gives an involution of $H$, this time fixing $H \backslash \left( C_1 \cup C_2\right)$ and partitioning the remainder into $C_1$ and $C_2$.  Furthermore, this construction can also be generalized to arbitrary subsets of isomorphic connected components by treating the sets of connected components as the sets of isomorphic vertices were treated above.

We now focus on more complex graph families that have been shown to be Sidorenko by methods other than sos.  We will show that since they are non-trivial squares, they could potentially admit sos proofs.  Furthermore, we offer square root constructions that could be used to start building such proofs and in some cases fully classify all square roots. 


\subsection{Square roots of graphs satisfying degree conditions}
We start with Theorem 1.1 from \cite{conlon2018sidorenko}, which states that Sidorenko's conjecture holds for bipartite graphs $H=(A\cup B, E)$ with the following property: $\binom{|A|}{r}\binom{r}{k}$ divides $d_k$ for each $1 \leq k \leq r$, where $d_k$ is the number of vertices of degree $k$ in $B$ and $r$ is the maximum degree of vertices in $B$. We show that all such graphs are non-trivial squares and give explicit constructions of square roots for these graphs. In fact, to be a non-trivial square it is sufficient that $d_k  > \binom{|A|}{k}$ for some $k$ or $d_1=\ldots=d_{r-1}=0$ and $d_r=\binom{|A|}{r}$.  These two conditions include all cases in \cite{conlon2018sidorenko}.

\begin{lemma}\label{lem1} Let $H=(A\cup B, E)$ be a bipartite graph where $d_k$ is the number of vertices of degree $k$ in $B$ and $r$ is the maximum degree of vertices in $B$.  If there exists $k\leq r$ such that $d_k>\binom{|A|}{k}$, then $H=[[G^2]]$ for some graph $G$ with $|V(H)|-1$ vertices that are all labeled and one unlabeled.
\end{lemma}

\begin{proof}
If there exists $k\leq r$ such that $d_k> \binom{|A|}{k}$, then since there are $\binom{|A|}{k}$ different potential neighborhoods for the $d_k$ vertices of degree $k$, this implies that at least two of those vertices must have the same neighborhood, say $v$ and $v'$. 

By Lemma \ref{twins}, $H$ is a non-trivial square and the square root construction in Lemma \ref{twins} may be used.
\end{proof}

\begin{lemma}\label{lem2}
Let $H=(A\cup B, E)$ be a bipartite graph where $d_k$ is the number of vertices of degree $k$ in $B$ and $r$ is the maximum degree of vertices in $B$. If $d_1=\ldots=d_{r-1}=0$ and $d_r=\binom{|A|}{r}$, then $H$ is a non-trivial square.
\end{lemma}

\begin{proof} Here we have a graph where $|B|=\binom{|A|}{r}$ and every vertex in $B$ has degree $r$. Note that if any two of them have the same neighborhood, then by Lemma \ref{twins}, we have that $H$ is a non-trivial square. So we may assume that each vertex of $B$ has a distinct neighborhood in $A$, in other words each of the $\binom{|A|}{r}$ subsets of $r$ vertices in $A$ is a neighborhood of exactly one vertex in $B$. 
Call such a graph an \emph{$(|A|,r)$ graph}. To lighten the notation, let $A=\{1,2,\ldots, m\}$, so $|A|=m$.

If $r=1$, then $|B|=|A|$ and every vertex in $B$ has degree 1, and no two have the same neighborhood. Thus, $H$ is a collection of disjoint edges. If $H$ contains more than one edge, by Lemma~\ref{twincc}, $H$ is a non-trivial square. Note that if $H$ is a single edge, then $H$ is a trivial square. However, Sidorenko's conjecture is trivially true in that case: $\uvedge-\uvedge\geq 0$. 

If $m=1$, then $r=0$ (and the graph consists of two isolated vertices, and is a non-trivial square by Lemma~\ref{twincc}) or $r=1$ (and the graph is a single edge, which is not a non-trivial square, but for which Sidorenko's conjecture is trivially true as seen above). 

If $m=2$, either $r=0$ (and the graph consists of three isolated vertices, and is a non-trivial square by Lemma~\ref{twincc}), $r=1$ (and the graph is two disjoint edges which is a non-trivial square by Lemma~\ref{twincc}) or $r=2$ (and the graph is a path of length 2 which is a non-trivial square as seen in the introduction).

If $m\geq 3$, we build a square root for $H$ as follows. Consider the following partially labeled bipartite graph $G=(\bar{A}\cup \bar{B}, \bar{E})$ in Figure \ref{fig:degcond}. Let $\bar{A}$ contain $m-2$ vertices labeled $a_3, \ldots, a_m$ and one unlabeled vertex which we will call $a^*$. Let $L_{\bar{A}}=\{a_3, \ldots, a_m\}$ be the set of labeled vertices in $\bar{A}$. Let $\bar{B}$ contain a set of $\binom{m-2}{r}$ labeled vertices, called $L_1$, an additional set of $\binom{m-2}{r-2}$ labeled vertices, called $L_2$, and a set of $\binom{m-2}{r-1}$ unlabeled vertices, called $U_{\bar{B}}$. The set of edges $\bar{E}$ is partitioned into four types. The first type is composed of edges between $L_{\bar{A}}$ and $L_1$ to form the graph $(m-2,r)$, that is, every vertex in $L_1$ is adjacent to a different set of $r$ vertices in $L_{\bar{A}}$. The second type is composed of $\binom{m-2}{r-1}$ edges between vertices in $L_2$ and $a^*$. The third type is composed of edges between $L_{\bar{A}}$ and $U_{\bar{B}}$ so to form a $(m-2, r-1)$ graph, i.e., every vertex in $U_{\bar{B}}$ is adjacent to a different set of $r-1$ vertices in $L_{\bar{A}}$. The fourth type are the $\binom{m-2}{r-1}$ edges between $U_{\bar{B}}$ and $a^*$. We claim that $[[G^2]]$ gives us an $(m, r)$ graph, i.e., $[[G^2]]=H$.

\begin{figure}[h]
    \centering
    \includegraphics[scale=0.33]{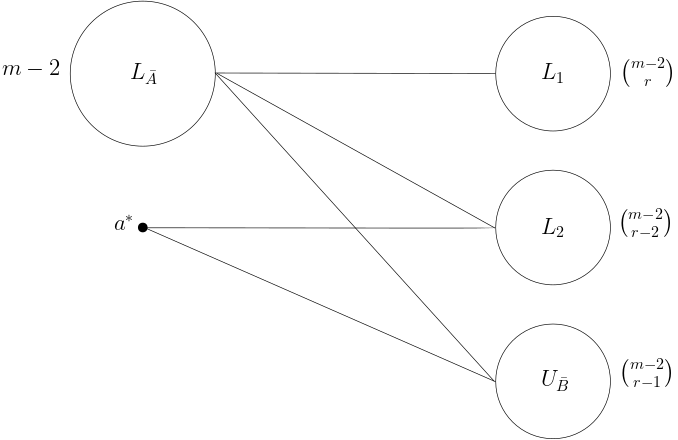} \quad \includegraphics[scale=0.33]{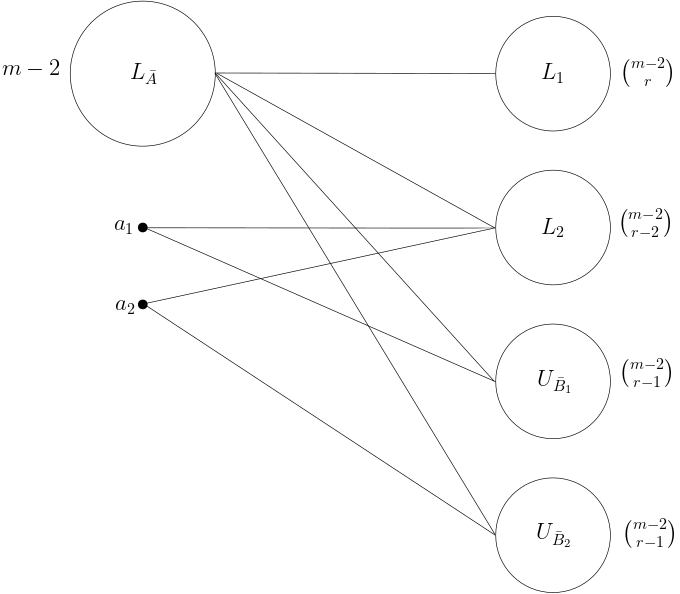}
    \caption{Cartoon of $G$ and $G^2$}
    \label{fig:degcond}
\end{figure}

First note that $G^2$ is still a partially labeled bipartite graph.   In the square, the $\bar{A}^2$ part is now the original $m-2$ labeled vertices with two unlabeled copies of $a^*$, which for ease of notation, we will call $a_1$ and $a_2$.  The $\bar{B}^2$ part is now the original sets of labeled vertices $L_1$ and $L_2$, and two additional copies of $U_{\bar{B}}$, which for ease of notation we will call $U_{\bar{B}_1}$ and $U_{\bar{B}_2}$.  This gives a total of $m$ vertices in the $\bar{A}^2$ part and $\binom{m-2}{r}+\binom{m-2}{r-2}+2\binom{m-2}{r-1}=\binom{m}{r}$ vertices in the $\bar{B}^2$ part.  

Furthermore, the edges are assigned correctly.  It is enough to show that each subset of $\bar{A}^2$ is the neighbor set of exactly one vertex in $\bar{B}^2$.  Any subset of $a_3, \ldots a_m$ is still the neighbor set of exactly one vertex in $L_1$ and no other vertices, since other vertices all neighbor at least one of $\{a_1, a_2\}$.  Any subset of $\{a_2, a_3, a_4, \ldots, a_m\}$ containing $a_2$ is a neighbor set of the vertex in $U_{\bar{B}_2}$ whose neighbors in $G$ were that same set with $a^*$ instead of $a_2$.  Similarly any subset of $\{a_1, a_3, a_4, \ldots, a_m\}$ containing $a_1$ is a neighbor set of the vertex in $U_{\bar{B}_1}$ whose neighbors in $G$ were that same set with $a^*$ instead of $a_1$.  Finally, any subset of $\{a_1, a_2, a_3, \ldots a_m\}$ containing both $a_1$ and $a_2$ is a neighbor set of the vertex in $L_2$ whose neighbors in $G$ were that same set with $a^*$ replacing both $a_1$ and $a_2$.
\end{proof}

These two lemmas cover every case in Theorem 1.1 of \cite{conlon2018sidorenko}:

\begin{corollary}
Let $H=(A\cup B, E)$ be a bipartite graph such that $\binom{|A|}{r}\binom{r}{k}$ divides $d_k$ for each $1 \leq k \leq r$ where $d_k$ is the number of vertices of degree $k$ in $B$ and $r$ is the maximum degree of vertices in $B$. Then $H$ is a non-trivial square.
\end{corollary}
\begin{proof}
For $k<r$, $\binom{|A|}{r}\binom{r}{k}$ divides $d_k$ implies either $d_k=0$ or $d_k\geq \binom{|A|}{r}\binom{r}{k}>\binom{|A|}{k}$.  If there exists $k$ for which the latter is true, Lemma \ref{lem1} applies.  If the former is true for all $k$, then since $d_r$ is non-zero by definition we either have $d_r=\binom{|A|}{r}\binom{r}{r}=\binom{|A|}{r}$ or else $d_r>\binom{|A|}{r}$.  If the first, then Lemma \ref{lem2} applies.  If the second, then Lemma \ref{lem1} holds.

\end{proof}


\subsection{Square roots of subdivisions of graphs}

For an unlabeled graph $H$, let $H^{\textup{div}}$ be its subdivision where every edge becomes a path of length two. Theorem 1.3 in \cite{conlonkimleelee} states that Sidorenko's conjecture holds for the subdivision of $K_k$ for $k\geq 2$. Furthermore, it states that if Sidorenko's conjecture holds for arbitrary $H$, then it also holds for $H^{\textup{div}}$.  (Note that \cite{conlonkimleelee} uses the notation $H_1$ for $H^{\textup{div}}$.) In this section, we show that the subdivision of $K_k$ is a non-trivial square, and that if $H$ is a non-trivial square, then $H^{\textup{div}}$ is a non-trivial square too.  Thus the Sidorenko results for these classes of graphs could potentially be proved using sos.

\begin{lemma}
$K_k^{\textup{div}}$ is a non-trivial square for all $d\geq 2$.
\end{lemma}

\begin{proof}
First observe that $K_2^{\textup{div}}=\uHtwo=\left[\left[ \left(\vedge{1}{}\right)^2\right]\right]$ and $K_3^{\textup{div}}$ is a 6-cycle, which is equal to $\left[\left[\left(\pthree{1}{}{}{2}\right)^2\right]\right]$. 

We now give a construction for when $k\geq 4$. For convenience, we call the vertices in the non-subdivided graphs the  \emph{original} vertices, and we call the remaining vertices of $K_k^{\textup{div}}$ the \emph{new} vertices. We construct a square root of $K_k^{\textup{div}}$ as follows.

Start by labeling completely $K_{k-2}^{\textup{div}}$.  Then add one unlabeled vertex $u$. Put a path of length two (i.e., a subdivided edge) between $u$ and each original vertex of $K_{k-2}^{\textup{div}}$. Furthermore, add a labeled vertex $\ell$ that is adjacent to $u$. Call this graph $G$. For example, here is the construction for $K_4^{\textup{div}}$ and $K_5^{\textup{div}}$:

$$\left[\left[\left(\Kfourdconstruction{1}{2}{3}{4}\right)^2\right]\right]=K_4^{\textup{div}} \textup{ and } \left[\left[\left(\Kfivedconstruction{7}{1}{6}{5}{2}{4}{3}\right)^2\right]\right]=K_5^{\textup{div}}$$

Looking at $G^2$, the labeled $K_{k-2}^{\textup{div}}$ remains, and we see that $u$ and its copy, say $u'$, have become the last two original vertices of $K_k^{\textup{div}}$ since there is a subdivided edge from every original vertex in $K_{k-2}^{\textup{div}}$ and $u$ and $u'$ respectively, and $\{u, l\}$ and $\{l, u'\}$ forms a subdivided edge between $u$ and $u'$.  Thus $[[G^2]]=K_k^{\textup{div}}$.
\end{proof}

 Before the next lemma, we make an observation about subdivisions of labeled graphs in general.  For a labeled graph $F$, let $F^{\textup{div}}$ be its subdivision where a vertex is added between any two vertices of $F$ to create a path of length two, as before. The difference here is that some of these new vertices will be labeled, namely the ones that are added between two labeled vertices of $F$. The new labeled vertices have a label that consists of the pair of labels of the two adjacent vertices.

It's easy to see that if $H_1=[[F_1^2]]$, then $H_1^{\textup{div}}=[[(F_1^{\textup{div}})^2]]$. Moreover, by labeling new vertices consistently, cross-products are also consistent with this idea: $[[F_1 F_2]]^{\textup{div}} = [[F_1^{\textup{div}} F_2^{\textup{div}}]]$.

\begin{example}
Here is an example of this observation:
$$H_1=\ucfour= \left[\left[\left(\Htwo{}{1}{2}\right)^2\right]\right] \textup{ and } H_1^{\textup{div}} = \uceight = \left[\left[\left(\Htwod{1}{2}\right)^2\right]\right] $$
\end{example}

We are now ready to state and prove the second subdivision lemma.

\begin{lemma}
If $H-\uvedge^{|E(H)|}$ is sos (and thus $H$ has Sidorenko's property), then $H^{\textup{div}} - \uvedge^{2|E(H)|}$ is also sos (and thus $H^{\textup{div}}$ has Sidorenko's property). 
\end{lemma} 

\begin{proof}
If $H- \uvedge^{|E(H)|}$ is sos, then without loss of generally we can assume that $H-\uvedge^{|E(H)|} = \sum_{i} \alpha_i [[ (F_{i1} - F_{i2})^2]]$ for some set of graphs $\{F_{i1}, F_{i2}\}$. Then, $H^{\textup{div}} - \uHtwo^{|E(H)|} = \sum_{i} [[(F_{i1}^{\textup{div}} - F_{i2}^{\textup{div}})^2]]$. Since $\uHtwo-\uvedge\uvedge=[[(\vedge{1}{} - \vedge{2}{})^2]]$,  this is enough to show that $F^{\textup{div}} - \uvedge^{2|E(H)|}$ is sos. 
\end{proof}
\subsection{Square roots of Cartesian products}

\begin{figure}[h]
     \centering
     \includegraphics[scale=0.3]{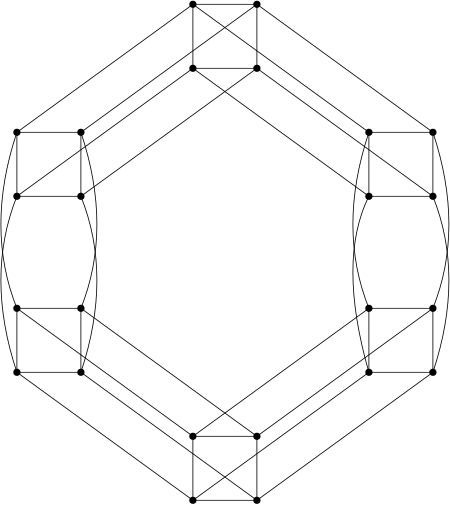}
     \caption{The Cartesian product $C_4 \square C_6$}
     \label{fig:cart}
 \end{figure}

 We now discuss the family of graphs described in Theorem 1.4 of \cite{conlonkimleelee}.  They consider the Cartesian product of an arbitrary Sidorenko graph $H$ with an arbitrary even cycle $C_{2k}$, denoted as $H \square C_{2k}$.  In other words, the vertex set is $\{(v, i) : v \in V(H), 1\leq i \leq 2k\}$ and there is an edge between $(u,i)$ and $(v,j)$ when $i=j$ and $u$ is adjacent to $v$ in $H$, or $u=v$ and $i-j \equiv 1 \mod{2k}$. (See Figure \ref{fig:cart} for an example.) Theorem 1.4 states that if $H$ is Sidorenko, then so is $H \square C_{2k}$.  Here we prove that any graph of the form $H \square C_{2k}$ is a non-trivial square (regardless of the properties of $H$).  
 
 We also fully classify square roots of these types of graphs.  If $H$ is a trivial square we show that $H \square C_{2k}$ has a unique square root.  If $H$ is a non-trivial square we show that all square roots of $H \square C_{2k}$ are constructed from square roots of $H$, apart from the construction used for trivial square $H$.

\begin{theorem}\label{default}
For any arbitrary (not necessarily Sidorenko) $H$, and any $k\geq 2$, $H \square C_{2k}$ is a non-trivial square.  A square root construction is given which holds for any arbitrary $H$ and $k\geq 2$.
\end{theorem}

\begin{theorem}\label{square}
If $[[A^2]]=H \square C_{2k}$, then $A$ is the square root construction from Theorem \ref{default} or else $A=B \square C_{2k}$ for some $B$ such that $[[B^2]]=H$.
\end{theorem}

In fact, Theorem \ref{square} has the following corollary:
\begin{corollary}
If $H$ is a trivial square,  the only possible square root of $H \square C_{2k}$ is the graph constructed in Theorem \ref{default}.
\end{corollary}

Taken together, these greatly restrict the possibilities for an sos proof of Theorem 1.4 in \cite{conlonkimleelee}.  Any such proof would have to rely on square roots for $H$, or the construction given in Theorem \ref{default}.  Expanding a Sidorenko graph $H$ to a Sidorenko graph $H \square C_{2k}$ in fact adds very little to the proof possibilities.

In the following discussion, we use some general notation for squaring graphs.  Suppose we are interested in some graph $A$ and its square $A^2$.  Consider the labeled and unlabeled vertices in $A$.  The labeled vertices in $A$ give rise to labeled vertices in $A^2$ while each unlabeled vertex in $A$ gives rise to two unlabeled \emph{twin vertices} in $A^2$.  In fact, the labeled vertices correspond to the vertices fixed by the involution of Lemma 2.8 of \cite{BRST} (discussed earlier).  The twin unlabeled vertices each belong to the two separate groups in the partition generated by the involution.  Observe that vertices in one part of the partition only neighbor other vertices within that part or vertices fixed by the involution.  We will call these two parts of the partition \emph{unlabeled odd} and \emph{unlabeled even}.  

If we want to discuss an odd vertex and even vertex that arise from the same unlabeled vertex in $A$, i.e. a pair of twin vertices in $A^2$, we will call them $u_{2j-1}$ and $u_{2j}$ for some convenient $j$.  In particular, then, any labeled vertex neighboring $u_{2j-1}$ must also neighbor $u_{2j}$ and vice-versa.  At times it will be useful to think of the unlabeled vertices $A^2$ as \emph{originals} and \emph{copies}; we will assume the odd unlabeled vertices are the originals.

We also standardize some notation for $H \square C_{2k}$.  Recall that $H \square C_{2k}$ is made up of $2k$ copies of $H$, which we call $H_1, H_2, \ldots H_{2k}$.  Let $\varphi_{j,k}$ be the canonical isomorphism between $H_j$ and $H_k$ that sends $(v, i)$ to $(v,j)$.  In this language, $H \square C_{2k}$ then has cross-edges between $v \in H_i$ and $\varphi_{i, i+1}(v)$ in $H_{i+1}$, and between $v \in H_i$ and $\varphi_{i, i-1}(v)$ in $H_{i-1}$.  The indexing is modulo $2k$, i.e., there is an edge between $v \in H_1$ and $\varphi_{1, 2k}(v)$.

\begin{figure}[h]
     \centering
     \includegraphics[scale=0.3]{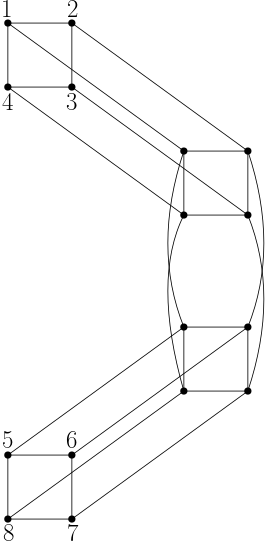}
     \caption{The square root construction of $C_4 \square C_6$ from Theorem \ref{default}}
     \label{fig:cartsr1}
 \end{figure}

With this notation in hand, we are ready to prove the two theorems.  First we prove Theorem \ref{default} with an explicit construction:

\begin{proof}[Proof of Theorem \ref{default}]
Start with $H \square C_{2k}$.  Label the vertices in $H_1$ and $H_{k+1}$, leave the vertices in $H_2$ through $H_k$ unlabeled, and remove the vertices in $H_{k+2}$ through $H_{2k}$. Squaring the remaining graph and unlabeling it yields $H\square C_{2k}$. (See Figure~\ref{fig:cartsr1} for an example.) \end{proof}

As Theorem \ref{square} is a negative result, its proof is more involved and uses the following lemma:
\begin{lemma}\label{identical}
If $[[A^2]]=H \square C_{2k}$, then either the sets of labeled, unlabeled odd, and unlabeled even vertices of $A$ are identical on each $H_i$, or else $A$ is the construction from Theorem \ref{default}.
\end{lemma}

\begin{proof}[Proof of Lemma \ref{identical}]

Suppose not, i.e., suppose we have an $A$ such that $[[A^2]]=H \square C_{2k}$ and the labelings of $H_i$ and $H_{i-1}$ are not identical for some $i$.  Without loss of generality, there is then a labeled vertex $\ell_a$ in $H_i$ such that $\varphi_{i, i-1}(\ell_a)$ is even unlabeled, say $u_2$.  In this case, its twin $u_1$ is also a neighbor of $\ell_a$.  Note that it's not possible that $u_1 \in H_{i-1}$ since $\ell_a$ already has a neighbor in $H_{i-1}$.  Therefore $u_1 \in H_i$ or $u_1 \in H_{i+1}$.  

First suppose $u_1$ in $H_i$.  Since $\ell_a$ neighbors $u_1$ in $H_i$, we know that $u_2$ must neighbor $\varphi_{i, i-1}(u_1)$ and by construction so must $u_1$.  The unlabeled vertices only neighbor other unlabeled vertices of the same parity; therefore $\varphi_{i, i-1}(u_1)$ must be a labeled vertex.  Call it $\ell_b$.

Now consider $\varphi_{i-1, i-2}(u_2)$.  This must be an unlabeled vertex, since if it were labeled it would also neighbor $u_1$ which is in $H_i$.  As discussed above, it must be even since it neighbors $u_2$.  Call $\varphi_{i-1, i-2}(u_2)=:u_4$.  Then we know that its twin, $u_3$, must neighbor $u_1$.  Therefore $u_3 \in H_{i+1}$ or $u_3 \in H_{i+2}$.  We also know that $u_4$  neighbors $\varphi_{i-1, i-2}(\ell_b)$.  If $\varphi_{i-1, i-2}(\ell_b)$ were labeled, it would have to neighbor $u_3$, which is impossible since they are in $H_{i-2}$ and $H_{i+1}$ or $H_{i+2}$.  Therefore $\varphi_{i-1, i-2}(\ell_b)$ is unlabeled, and must be even again, call it $u_6$.  Again its twin, $u_5$, must neighbor $u_3$.  Since $\ell_b$ neighbors $u_6$, it must also neighbor $u_5$.  The only way for this to happen, since $\ell_b \in H_{i-1}$ and $\ell_b$ already neighbors $u_1 \in H_i$, is that $u_3 \in H_i$ and $u_5=\varphi_{i, i-1}(u_3)$.  This subgraph is illustrated in Figure \ref{fig:arg1}.
\begin{figure}[h]
    \centering
    \includegraphics[scale=.45]{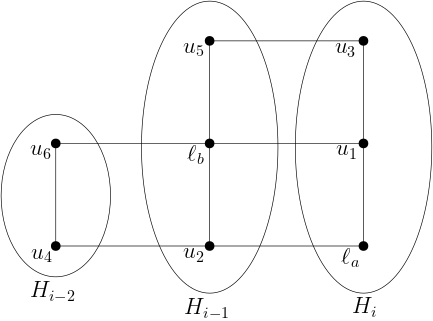}
    \caption{Configuration of vertices in $H_i, H_{i-1}, H_{i-2}$}
    \label{fig:arg1}
\end{figure}

Note that $u_5$, $\ell_b$, and $u_6$ are exactly in the same configuration as $u_1$, $\ell_a$, and $u_2$ that we originally began with.  Therefore the same argument can be repeated to generate a larger subgraph as in Figure \ref{fig:arg2}

\begin{figure}[h]
    \centering
    \includegraphics[scale=.45]{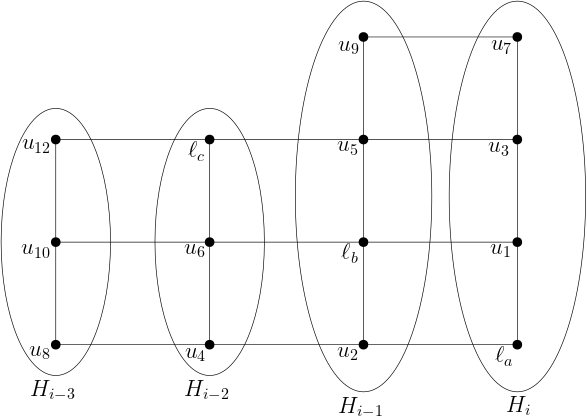}
    \caption{Configuration of vertices in $H_i, H_{i-1}, H_{i-2}, H_{i-3}$}
    \label{fig:arg2}
\end{figure}

In fact, repeating this argument $j$ times for arbitrary $j$ proves that  $\varphi_{i, i-j}(u_1)$ is an unlabeled even vertex.  In particular, letting $j=2k-1$ shows that $\varphi_{i, i-j}(u_1)=\varphi_{i, i+1}(u_1)$ must be an unlabeled even vertex.  On the other hand, it also must neighbor $u_1$, which is an unlabeled odd vertex.  This is a contradiction, so our original configuration of $\ell_a \in H_i$, $u_1 \in H_i$, $u_2 \in H_{i-1}$ is impossible.

The only remaining case is $\ell_a \in H_i$, $u_1 \in H_{i+1}$, $u_2 \in H_{i-1}$.  Since $H$ is a connected graph, $u_1$ must have a neighbor in $H_{i+1}$.  Note that this vertex cannot be labeled, as then it would have to neighbor $u_2 \in H_{i-1}$ as well.  Therefore it is odd unlabeled, for convenience, call it $u_3$.  Then $\varphi_{i+1, i-1} (u_3)$ must neighbor $u_2$ since $\varphi_{i+1, i-1} (u_1)=u_2$.  Similarly, $\varphi_{i+1, i}(u_3)=\varphi_{i-1, i}(u_4)$ must neighbor $\varphi_{i+1, i} (u_1)=\ell_a$.  Since this vertex neighbors both odd and even vertices, it must be labeled.  Call it $\ell_b$.

We have now shown that any vertex in the neighborhood of $u_1$ in $H_{i+1}$ must be unlabeled odd, and by symmetry any vertex in the neighborhood of $u_2$ in $H_{i-1}$ must be unlabeled even, and any vertex in the neighborhood of $\ell_a$ in $H_i$ must be labeled.  Furthermore these vertices in $H_{i+1}$ must be twins of these vertices in $H_{i-1}$.  Repeating this argument $D$ times where $D$ is the diameter of $H$ shows that in fact all of $H_i$ is labeled, all of $H_{i+1}$ is unlabeled odd, and all of $H_{i-1}$ is unlabeled even.  Furthermore, $H_{i+1}$ is the twin of $H_{i-1}$.

Now consider $H_{i+2}$.  If any vertex $v$ in $H_{i+2}$ were labeled, it would have to neighbor both $\varphi_{i+2, i+1}(v)$, which is in $H_{i+1}$, and that vertex's twin $\varphi_{i+1, i-1}\left(\varphi_{i+2, i+1}(v)\right)$ in $H_{i-1}$.  This is only feasible if $2k=4$, in which case we have the construction from Theorem \ref{default}.

Otherwise $H_{i+2}$ must be fully unlabeled odd, and similarly $H_{i-2}$ must be fully unlabeled even.  Furthermore, since $H_{i+1}$ and $H_{i-1}$ are twins, $H_{i+2}$ and $H_{i-2}$ must also be twins. Repeating this argument $k-1$ times shows that in fact we have the construction from Theorem \ref{default}.
\end{proof}

\begin{figure}[h]
     \centering
     \includegraphics[scale=0.3]{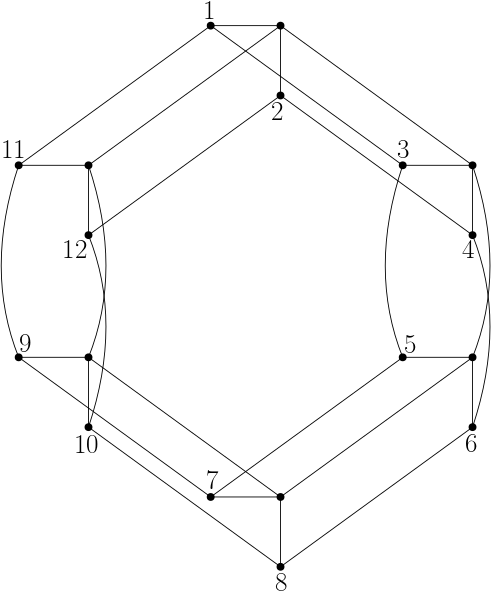}
     \caption{The square root construction of $C_4 \square C_6$ from Theorem \ref{square}}
     \label{fig:cartsr2}
 \end{figure}

With Lemma \ref{identical} in hand, we are now ready to prove Theorem \ref{square}.  
\begin{proof}[Proof of Theorem \ref{square}]
Suppose that some $A$ such that $[[A^2]]=H \square C_{2k}$ is in the first category, i.e., where each copy of $H$ has identical labelings.  Consider $A$ restricted to the vertices in $H_1$.  Recall our convention that the odd unlabeled vertices in $H$ are the originals from $A$ while the even unlabeled vertices are the \emph{copies}. Thus, $A$ restricted to the vertices in $H_1$ consists of the labeled and odd unlabeled vertices in $H_1$.  Call this subgraph $A_1$.

Now consider $[[A_1^2]]$.  (For an image of the following argument, see Figure \ref{fig:arg3}.)  Suppose there is some vertex outside of $H_1$ in $H \square C_{2k}$, without loss of generality, say in $H_2$.  Any such vertex must be even by construction.  Choose a vertex in $[[A_1^2]] \cap H_2$ that neighbors a vertex in $H_1$.  By the above and the fact that all $H_i$ are identically labeled, we already know that this vertex and its neighbor in $H_1$ must both be even.  Call the even vertex in $[[A_1^2]] \cap H_2$ $u_4$, and its neighbor in $H_1$ $u_2$.

Now consider the shortest path between $u_2$ and $u_1$ in $H_1$.  By construction, we know it must consist of even vertices $u_{2i_1}, u_{2i_2}, \ldots u_{2i_d}$, then a labeled vertex $\ell$, then the odd vertices $u_{2i_1-1}, u_{2i_2-1}, \ldots u_{2i_d-1}$, then $u_1$.  In other words it arises from the shortest path between $u_1$ and any labeled vertex in $A_1$.  

Look at what $\varphi_{1,2}$ does to this path.  We know that $\varphi_{1,2}(u_2)=u_4$ and  $$\varphi_{1,2}(u_{2i_1}), \varphi_{1,2}(u_{2i_2}), \ldots \varphi_{1,2}(u_{2i_d})$$ is a path of even vertices between $u_4$ and $\varphi_{1,2}(\ell)$.  The twins of these vertices must form a path of odd vertices between $\varphi_{1,2}(\ell)$ and $u_3$.  On the other hand, $\varphi_{1,2}$ maps $u_3$ to some odd vertex in $H_2$, call it $u_5$; this is the penultimate vertex on the twin path from $\varphi_{1,2}(\ell)$ and $u_3$, and the final vertex on the path contained in $H_2$.  This gives an even-vertex path of length $d$ between $\varphi_{1,2}(\ell)$ and $\varphi_{1,2}(u_2)$ but an odd-vertex path of length $d-1$ between $\varphi_{1,2}(\ell)$ and $\varphi_{1,2}(u_3)$.  Therefore, there must be an odd-vertex path of length $d-1$ between $\ell$ and $u_3$.  But then there must be an even-vertex path of length $d-1$ between $\ell$ and $u_4$, made up of the twins of the odd-vertex path!  By construction, the shortest path between $\ell$ and $u_2$ has length $d$ so the path of length $d-1$ between $\ell$ and $u_4$ must not contain $u_2$.  However, since $\ell$ is in $H_1$ and $u_4$ is in $H_2$, any path between them must contain $u_2$, which is a contradiction.

\begin{figure}
    \centering
\includegraphics[scale=.4]{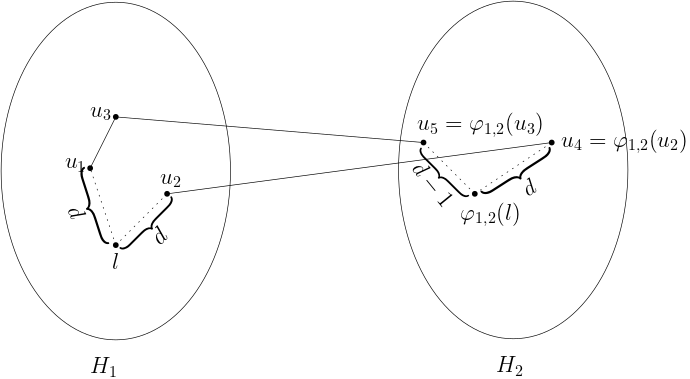}
     \caption{Subgraph from Proof of Theorem \ref{square}}
    \label{fig:arg3}
\end{figure}
This proves that $[[A_i^2]] \subseteq H_i$ for all $i$.  Therefore $[[A_i^2]]=H_i$; if it were true that $[[A_i^2]]$ were properly contained in $H_i$ for some $i$, then there would have to be some $j$ such that $[[A_j^2]] \cap H_i$ is nonempty, which contradicts the above.  We have now shown that $[[A_i^2]]=H_i$ for all $i$, so in fact $A=B \square C_{2k}$ for $B$ such that $[[B^2]]=H$, and the theorem is proved. (See Figure~\ref{fig:cartsr2} for an example of the construction in this theorem.)
\end{proof}

\section{Computational results}\label{sec:computationalresults}

We now provide some computational results to classify cases on at most 7 edges for which there is an sos proof for Sidorenko's conjecture. We will see that trivial squares are not the only obstacle to sos certificates.

Let $\mathcal{P}_k$ denote the poset of all $k$-edge bipartite graphs ordered by the relation $\succeq$  where $H_1\succeq H_2$ if $H_1-H_2 \geq 0$ is valid when evaluated on all target graphs $G$, i.e., $t(H_1;G)-t(H_2;G)\geq 0$ for all graphs $G$, and can be certified with an sos.  We explain in Section \ref{subsec:coverrelations} how to compute such posets. Figures \ref{fig:P5} and \ref{fig:P6}  give the Hasse diagrams of $\mathcal{P}_5$ and $\mathcal{P}_6$. The Hasse diagram for $\mathcal{P}_7$ is too big to be printed, and so can be viewed online instead: \url{https://people.math.umass.edu/~raymond/poset7whole.pdf}. In Figures \ref{fig:P7comp} and \ref{fig:P7incomp}, we give the Hasse diagrams of the subposets of $\mathcal{P}_7$ containing elements comparable and incomparable with $\uvedge^7$ respectively.

\begin{center}
\begin{figure}
    \centering
    \includegraphics[scale=0.5]{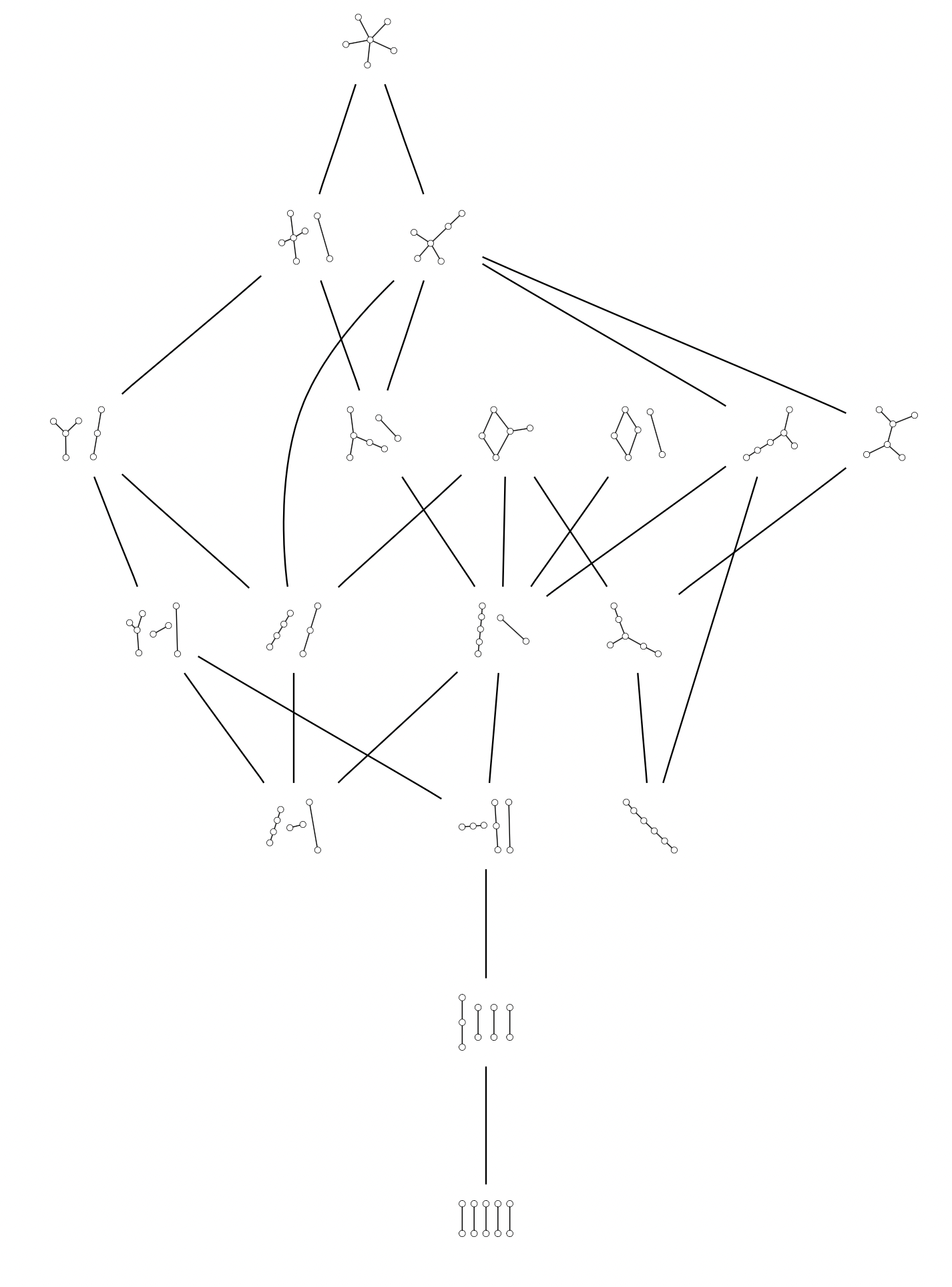}
    \caption{Hasse diagram of $\mathcal{P}_5$}
    \label{fig:P5}
\end{figure}

\begin{figure}
    \centering
    \includegraphics[scale=0.82]{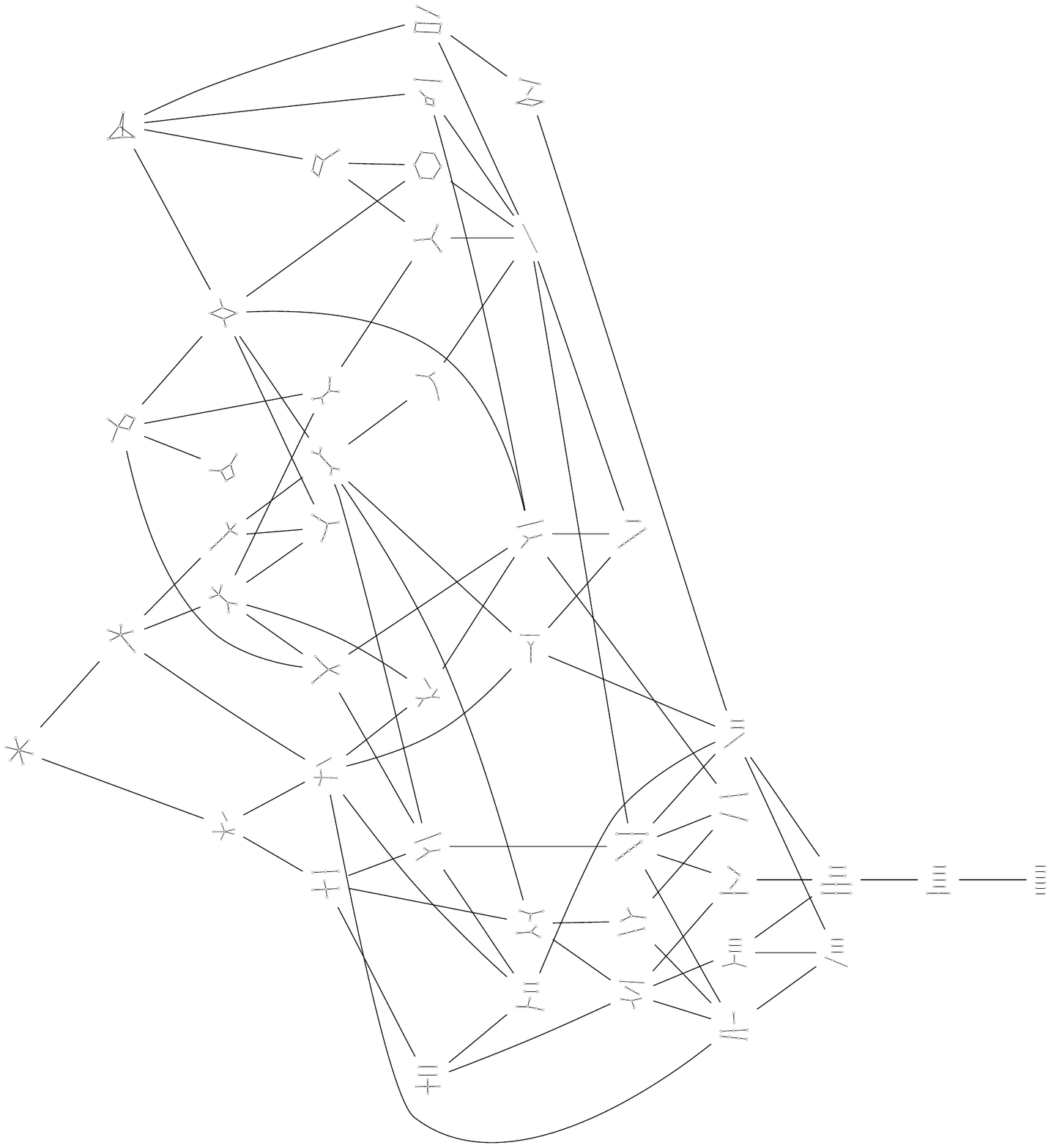}
    \caption{Hasse diagram of $\mathcal{P}_6$}
    \label{fig:P6}
\end{figure}

\begin{figure}
    \centering
    \includegraphics[scale=0.3]{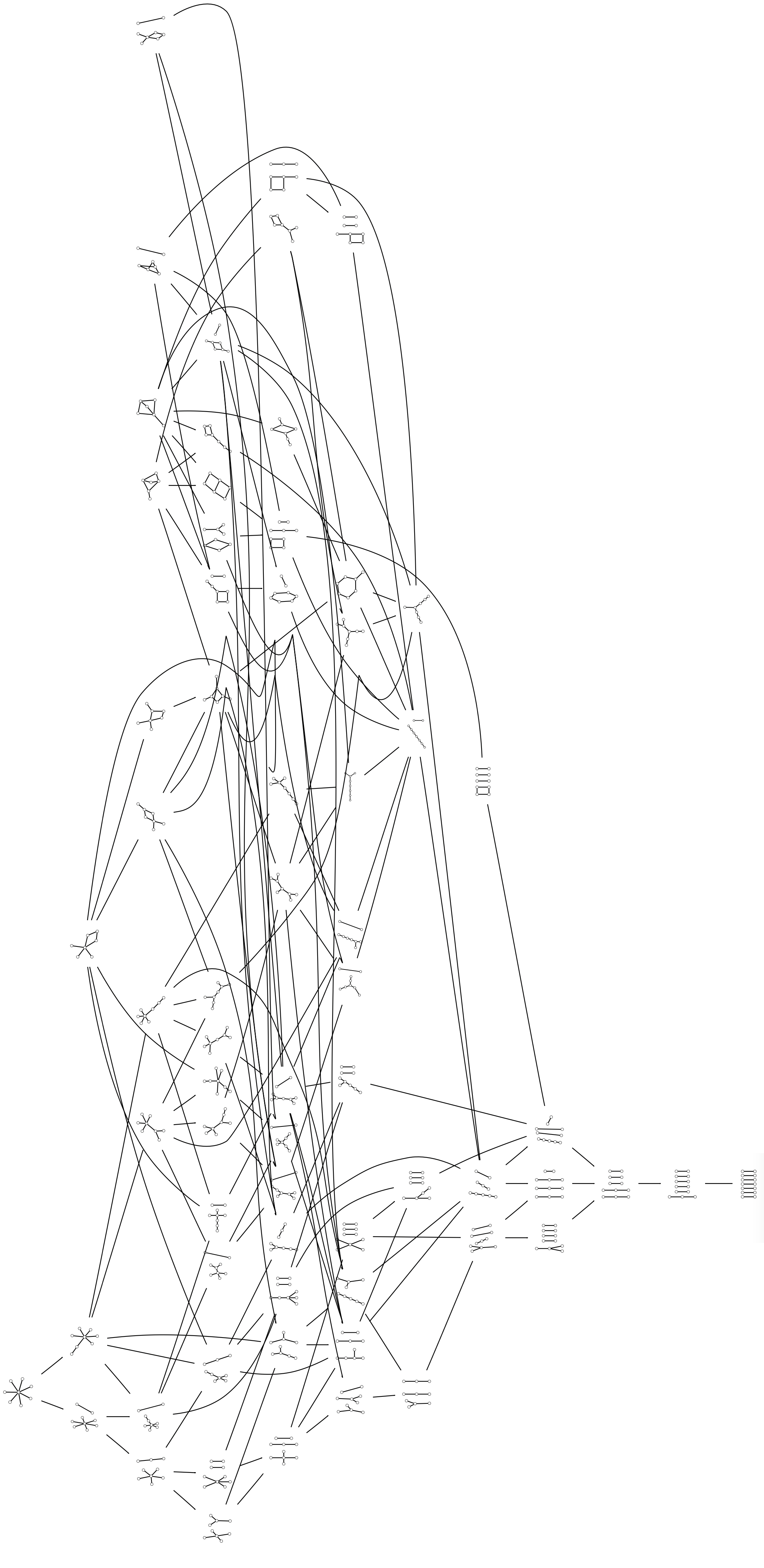}
    \caption{Hasse diagram of the subposet of $\mathcal{P}_7$ containing graphs that can be shown to be Sidorenko with an sos}
    \label{fig:P7comp}
\end{figure}

\begin{figure}
    \centering
    \includegraphics[scale=0.45]{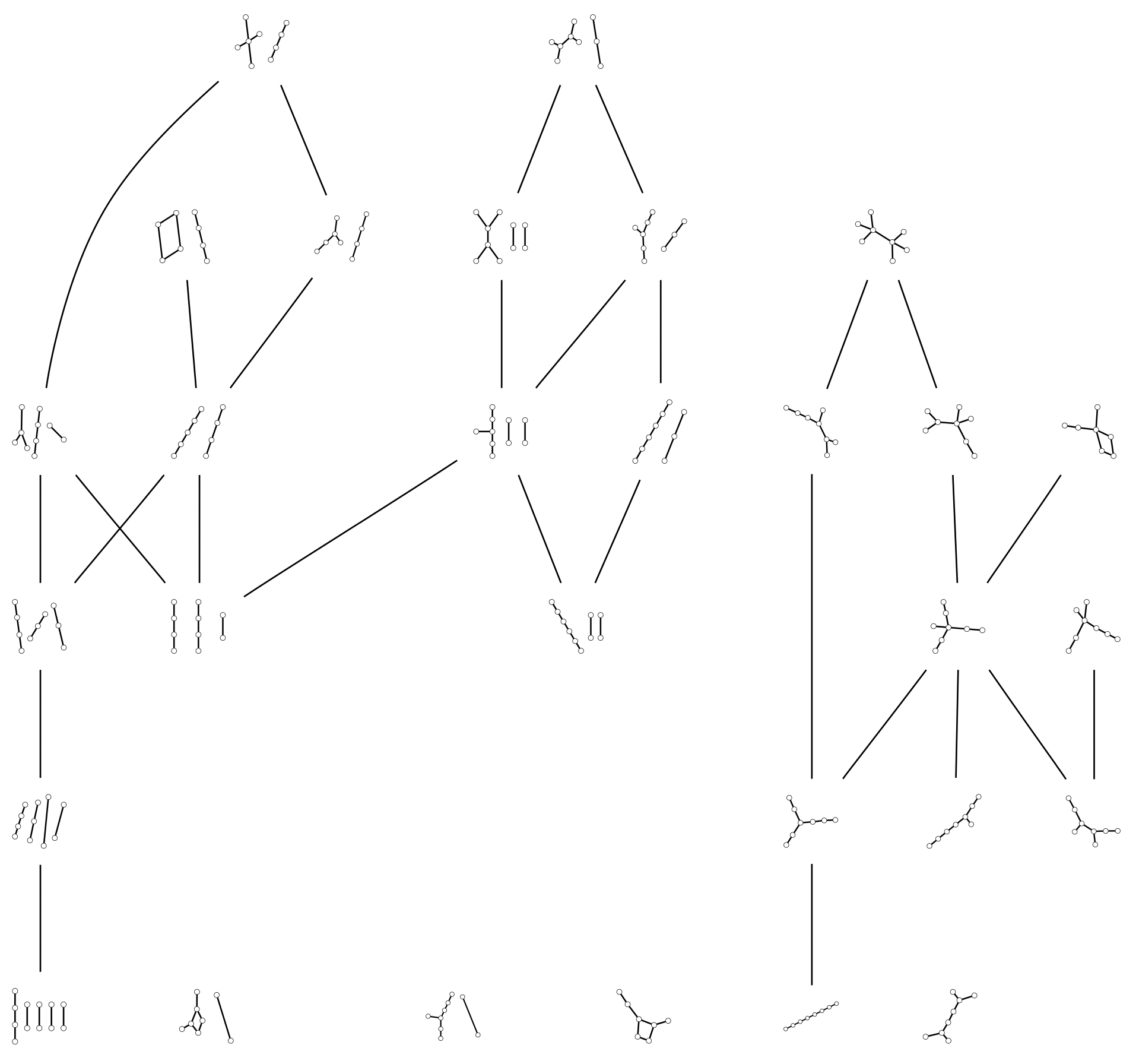}
    \caption{Hasse diagram of the subposet of $\mathcal{P}_7$ containing graphs that cannot be shown to be Sidorenko with an sos}
    \label{fig:P7incomp}
\end{figure}
\end{center}

Observe that if $H_1\succeq H_2$ and $H_2\succeq H_3$, then $H_1\succeq H_3$. Indeed, if $H_1\succeq H_2$ and $H_2\succeq H_3$, then $H_1-H_2=\sum_{i}[[a_i^2]]$ and $H_2-H_3=\sum_{j} [[b_j^2]]$ where the $a_i$'s and $b_j$'s are graph combinations, which implies that $H_1-H_3=(H_1-H_2)+(H_2-H_3)=\sum_{i}[[a_i^2]] + \sum_j [[ b_j^2]]$, and so $H_1 \succeq H_3$.

Therefore, a graph $H$ with $k$ edges is such that $H-\uvedge^{|E(H)|}$ has an sos proof if and only if there is a chain from $\uvedge^k$ to $H$ in $\mathcal{P}_k$.

Let $\mathcal{Q}_k$ denote the poset of all $k$-edge bipartite graphs ordered by the relation $\rhd$ where $H_1 \rhd H_2$ if $H_1-H_2\geq 0$ is valid when evaluated on all target graphs $G$. Note that $H_1\succeq H_2$ implies that $H_1 \rhd H_2$, i.e., $\succeq$ is a subrelation of $\rhd$, but the converse isn't true. Sidorenko's conjecture claims that $\uvedge^k$ is the unique minimal element of $\mathcal{Q}_k$.

In $\mathcal{P}_5$, we see that there are five bipartite graphs $H$ with five edges for which $H-\uvedge^{5}$ cannot be written as a sum of squares: $\upfive$, $\upthree\ \uvedge \ \uvedge$, $\upthree \ \uHtwo$,  $\uchristmastree$ and $\uweights$. From \cite{BRST}, it was already known that $\upfive - \uvedge^5$ cannot be written as a sum of squares since $\upfive$ is a trivial square. It is also unsurprising that $\upthree\ \uvedge \ \uvedge - \uvedge^5$ and $\upthree \ \uHtwo -\uvedge^5$ cannot be written as sums of squares since $\upthree$ is also a trivial square and so it is known that $\upthree -\uvedge^3$ has no sum of squares certificate. On the other hand, $\uchristmastree$ and $\uweights$ are not trivial squares since $$\uchristmastree=\left[\left[\left(\pthree{1}{2}{}{}\right)^2\right]\right]$$ and $$\uweights=\left[\left[\left(\pthree{}{1}{2}{}\right)^2\right]\right]=\left[\left[\left(\broom{}{1}{2}{3}{4}\right)^2\right]\right].$$

However, note that there is a chain from $\upfive$ to $\uchristmastree$ and $\uweights$, that is, we have sos certificates for $\uchristmastree-\upfive$ and $\uweights-\upfive$. This means that if there is some other type of certificate of nonnegativity for $\upfive-\uvedge^5$ (which in this case there is, see for example \cite{blakleyroy}, \cite{london} and \cite{mulhollandsmith}), then we still have a way of showing the nonnegativity of  $\uchristmastree-\uvedge^5$ and $\uweights-\uvedge^5$ by combining that other type of certificate with sums of squares.

\textbf{Problem:} Characterize the minimal elements of $\mathcal{P}_k$. 

Having a characterization of the minimal elements of $\mathcal{P}_k$ would make clear what classes of graphs other proof techniques should focus on.  We give the minimal elements for $k=6$ and $k=7$.

In $\mathcal{P}_6$, besides $\uvedge^6$, the minimal elements are $\upthree \uvedge \uvedge \uvedge$, $\upthree \ \upthree$, $\upfive \uvedge$, $\usquareteletubbytwoantennas$, $\ustaronetwothree$ which all contain at least one connected component that is a trivial square. 

In $\mathcal{P}_7$, besides $\uvedge^7$, the minimal elements are $$\upthree \uvedge \uvedge \uvedge \uvedge, \upthree \ \upthree \uvedge, \upfive \uvedge \uvedge, \usquareteletubbytwoantennas \uvedge, \ustaronetwothree \uvedge,$$ $$\upseven, \uspider, \uyseven,  \uuneventeletubby, \uloongweights.$$ Note that \uloongweights is not a trivial square, but all other graphs contain at least one connected component that is a trivial square. This tells us that the characterization of minimal elements of $\mathcal{P}_k$ will not simply rely on trivial squares.

\subsection{Sum of squares certificates of cover relations}\label{subsec:coverrelations}

To compute the posets in the previous section, we had to check whether $H_1-H_2$ had an sos certificate for every pair of bipartite graphs $H_1$, $H_2$ with $k$ edges. This could have been done with a semidefinite program, but we instead chose to use Proposition 2.10 from \cite{BRST} (rephrased below) to have a linear program.

\begin{proposition}[Blekherman, Raymond, Singh, Thomas, 2020] \label{prop:oldpaper}
A graph combination $H_1-H_2$ where $H_1$ and $H_2$ have the same number of edges can be written as an sos if and only if it has an sos certificate of the form $H_1-H_2=\sum \lambda_{ij}[[(F_i - F_j)^2]]$  where $\lambda_{ij} \geq 0$ and where $|E([[F_i^2]])|=|E([[F_j^2]])|=|E([[F_iF_j]])|$ for each $i,j$, i.e., where $F_i$ and $F_j$ have the same fully labeled edges and the same number of edges.
\end{proposition}

We first produce all partially labeled graphs whose squares have $k$ edges. We then produce their squares as well as the cross-products of any pair of partially labeled graphs that have the same fully labeled edges: these are the pairs corresponding to possible $F_i$, $F_j$. We can then create a linear program where the variables are the $\lambda_{ij}$ and where we simply check whether there exist nonnegative $\lambda_{ij}$'s that yield $H_1-H_2$. To do so, we have a constraint for every graph $H$ with $k$ edges where we add $\lambda_{ij}$ if $[[F_i^2]]=H$, $\lambda_{ij}$ if $[[F_j^2]]=H$, and $-2\lambda_{ij}$ if $[[F_iF_j]]=H$. We then set this expression to be equal to $1$ if $H=H_1$, $-1$ if $H=H_2$ and 0 otherwise.

In the Appendix are the sos certificates for the cover relations of the posets of the previous section. To get an sos proof for some $H-\uvedge^{|E(H)|}$, add the sos certificates for all cover relations on a chain from $\uvedge^{|E(H)|}$ and $H$ (again, if there is no such chain, then $H-\uvedge^{|E(H)|}$ cannot be written as an sos). In the next section, we discuss further why sums of squares fail in certain instances beyond trivial squares.

In the tables in the Appendix, we give sos certificates for $A-B$ where $A$ and $B$ are both bipartite and with 5, 6 or 7 edges, and where $A$ is connected. From the posets, we see that an sos for any case where $A$ and $B$ are both not connected graphs can be retrieved through the sos certificates of some smaller graphs. The first column contains $A$, and the second column contains $B$.  A particular $A,B$ might have more than one rows. Each row corresponds to some $\alpha_i[[(H_{i1}-H_{i2})^2]]$ where the third column contains $\alpha_i$, the fourth $H_{i1}$ and the fifth $H_{i2}$, and so that $A-B=\sum_i \alpha_i[[(H_{i1}-H_{i2})^2]]$ where we are summing over the rows corresponding to a particular $A,B$ in the table.

\section{Obstacles to sums of squares certificates}\label{sec:nosos}

From our computations, we see that certain inequalities cannot be certified with sos. We now give some ideas as to why that is the case. We already know from \cite{BRST} that Sidorenko's inequality for trivial squares cannot be written as an sos. However, as seen from the posets, there are other obstacles to sos proofs as we found non-trivial squares for which there was no sos proof. 

First recall from Proposition \ref{prop:oldpaper} that we only need to consider sums of binomials squared where the two graphs in any given binomial have exactly the same fully labeled edges and the same number of edges.

\begin{lemma}\label{lem:neg}
Let $A$ and $B$ be bipartite graphs with the same number of edges such that $A\geq B$. Consider a square root of $A$, say $H$, and consider a graph $F$ that has the same fully labeled edges as $H$ and the same number of edges as $H$. If there exists a graph $G$ such that  $(1-\alpha)t(A;G)+2\alpha t([[HF]];G)-\alpha t([[F^2]];G)- t(B;G)<0$, then there exists no sos $q$ and $\alpha > 0$ such that $A-B=\alpha[[(H-F)^2]]+q$.
\end{lemma}

\begin{proof}
Suppose not. Then $A-B=\alpha[[(H-F)^2]]+q$ which implies that $q=(1-\alpha)A+2\alpha[[HF]]-\alpha[[F^2]]-B$. Since $q$ is sos, this means that  $(1-\alpha^2)A+2\alpha[[HF]]-\alpha^2[[F^2]]-B\geq 0$ no matter on what graph we evaluate this expression. But we know that there exists a graph $G$ such that  $(1-\alpha)t(A;G)+2\alpha t([[HF]];G)-\alpha t([[F^2]];G)- t(B;G)<0$, so this is not the case.
\end{proof}

We now consider a specific case of the previous lemma using the partition(s) of the vertex sets of the bipartite graphs we are considering.

\begin{definition}
For a bipartite graph $G$, let $b(G)=\max|\{V_1(G) : V(G)=V_1(G)\cup V_2(G) \textup{ and } E(G)\in V_1(G)\times V_2(G)\}|-|V(G)|$. 
\end{definition}

Observe that $t(G;S_n)=O(n^{b(G)})$ as $n\rightarrow \infty$ where $S_n$ is the star with $n$ branches. 

\begin{lemma}\label{lem:star}
Let $A$ and $B$ be bipartite graphs with the same number of edges such that $A\geq B$. Consider a square root of $A$, say $H$, and consider a graph $F$ that has the same fully labeled edges as $H$ and the same number of edges as $H$. If $\max\{b([[HF]]), b(A)\} < b([[F^2]])$, then there exists no sos $q$ and $\alpha > 0$ such that $A-B=\alpha[[(H-F)^2]]+q$.
\end{lemma}

\begin{proof}
From Lemma \ref{lem:neg}, let's see what happens to $(1-\alpha^2)A+2\alpha[[HF]]-\alpha^2[[F^2]]-B\geq 0$ if we evaluate on $S_n$ as $n\rightarrow \infty$: we obtain $(1-\alpha)O(n^{b(A)})+2\alpha O(n^{b([HF]])})-\alpha O(n^{b([[F^2]])}) - n^{O(b(B))}$.  Since $\max\{b([[HF]]), b(A)\} < b([[F^2]])$, for $n$ large enough, we see that this expression is negative.  Thus no sos certificate can contain $\alpha[[(H-F)^2]]$. 
\end{proof}

We now add another helpful lemma.

\begin{lemma}\label{lem:crosstriv}
Let $A$ and $B$ be graphs with the same number of edges such that $A\geq B$. Consider the set of all square roots of $A$, say $\{H_i\}_i$, and for each $i$ consider the set of all graphs $\{F_{ij}\}_j$ that have the same fully labeled edges as $H_i$ and the same number of edges as $H_i$. Let $\mathcal{S}=\{(i,j)| [[H_i F_{ij}]] \textup{ is not a trivial square}\}$. Then there exists no sos $q$ and $\alpha_{ij} \geq0$ for $(i,j)\notin \mathcal{S}$ and $\alpha_{ij}=0$ for $(i,j)\in \mathcal{S}$ such that $A-B=\sum_{i,j} \alpha_{ij}[[(H_i-F_{ij})^2]]+q$.
\end{lemma}

\begin{proof}
Suppose not. Then $A-B=\sum_{i,j} \alpha_{ij}[[(H_i-F_{ij})^2]]+q$ where $q$ is sos, $\alpha_{ij}=0$ for all $(i,j)\in \mathcal{S}$, and $\alpha_{ij}\geq 0$ for all $(i,j)\not\in \mathcal{S}$. Since the set $\{H_i\}_i$ contains all square roots of $A$, by Lemma 2.4 of \cite{BRST}, we know that $A$ doesn't appear as a square or as a cross-product in $q$. So the previous equation can be rewritten as $q=-B-\sum_{i,j} \alpha_{ij} [[F_{ij}^2]] +\sum_{i,j} 2\alpha_{ij}[[H_iF_{ij}]]$, that is, we see that the only positive terms in the right-hand side have the form $[[H_iF_{ij}]]$ and are trivial squares. Therefore, by Theorem 1.2 of \cite{BRST}, $-B-\sum_{i,j} \alpha_{ij} [[F_{ij}^2]] +\sum_{i,j} 2\alpha_{ij}[[H_iF_{ij}]]$ cannot be an sos, a contradiction.  
\end{proof}

\begin{example}
Consider $\uchristmastree-\uvedge^5$. We show how the previous ideas restrict the existence of an sos. First observe that $\uchristmastree$ only has, up to labeling, a single square root, namely $\uchristmastree=[[ (\pthree{1}{2}{}{})^2]]$. Also note that $b(\uchristmastree)=-3$. We now list all graphs $F$ on three edges that have the same fully labeled edge as $\pthree{1}{2}{}{}$ and compute $b([[F]]^2)$ and $b([[F\pthree{1}{2}{}{}]])$.

$$\begin{array}{|c|c|c|c|c|}
\hline
F & [[F\pthree{1}{2}{}{}]] & b([[F\pthree{1}{2}{}{}]]) & [[F^2]]& b([[F^2]])\\
\hline
\pthree{2}{1}{}{} & \upfive & -3 & \uchristmastree & -3\\
\threestar{2}{1}{}{} & \uthreebroom & \textcolor{red}{-2} & \ufivestar & \textcolor{red}{-1}\\
\threestar{1}{2}{}{}  & \ulongbroom & \textcolor{red}{-2} & \ufivestar & \textcolor{red}{-1}\\
\pthree{}{1}{2}{} & \uchristmastree & -3 & \uweights & -3\\
\Htwo{1}{2}{} \vedge{3}{} & \upfour \uvedge &  \textcolor{red}{-3}  & \uthreestar \uHtwo & \textcolor{red}{-2} \\
\Htwo{2}{1}{} \vedge{3}{} & \ubroom \uvedge & \textcolor{red}{-3} & \uthreestar \uHtwo & \textcolor{red}{-2}\\
\Htwo{1}{2}{} \uvedge & \upfour \uvedge &  -3  & \uthreestar \uvedge \uvedge & -3 \\
\Htwo{2}{1}{} \uvedge & \ubroom \uvedge & -3 & \uthreestar \uvedge \uvedge & -3\\
\vedge{1}{2} \Htwo{}{3}{4} & \upthree \uHtwo & -3 & \usquare \uvedge & -3\\
\vedge{1}{2} \Htwo{}{3}{} & \upthree \uHtwo & -3 & \upfour \uvedge & -3\\
\vedge{1}{2} \Htwo{3}{}{} & \upthree \uHtwo & \textcolor{red}{-3} & \upfour \ufourstar & \textcolor{red}{-2}\\
\vedge{1}{2}\uHtwo & \upthree\uHtwo & -3 & \uHtwo\uHtwo \uvedge & -3\\
\vedge{1}{2} \vedge{3}{} \vedge{4}{} & \upthree \uvedge\uvedge & -4& \uHtwo \uHtwo \uvedge & -3 \\
\vedge{1}{2} \vedge{3}{} \vedge{}{} & \upthree \uvedge\uvedge & -4 & \uHtwo \uvedge\uvedge\uvedge & -4\\
\vedge{1}{2} \vedge{}{} \vedge{}{} & \upthree \uvedge\uvedge  & -4 & \uvedge\uvedge\uvedge\uvedge\uvedge & -5\\

\hline
\end{array}$$

By Lemma \ref{lem:star}, we can remove some $F$'s, namely the ones for which the third and fifth column are in red. Further, from Lemma \ref{lem:crosstriv}, we know that $[[(\pthree{1}{2}{}{}-\pthree{2}{1}{}{})^2]]$ cannot be the only binomial remaining present as its cross-product is a trivial square. These considerations restrict the search space for an sos proof. 

\end{example}


\section{Conclusions and Future Work}
Using ideas from past work on gluing algebras as a jumping-off point, we have analyzed the prospects of sos proofs for Sidorenko's conjecture in many settings.  We gave a few general categories of graphs that are guaranteed to have square roots.  We also showed that many graphs already known (by other methods) to be Sidorenko are non-trivial squares, and therefore may admit an sos proof.  Furthermore, we described square roots for these graphs, which are the necessary building blocks for such a proof.

Complementing this work, we gave some computational results. Using a linear program, for $m\leq 7$, we computed a poset containing all bipartite graphs with $m$ edges that shows exactly which bipartite graphs $H_1$, $H_2$ on $m$ edges are such that $H_1-H_2$ has an sos certificate. These computations yielded many examples where there is no sos certificate even though $H_1$ is not a trivial square. We believe that better understanding this poset would be beneficial to making progress on Sidorenko's conjecture. In particular, characterizing what the minimal elements of this poset are would make clear for what classes of graphs other techniques will be needed to solve the conjecture.
We also discussed theoretical barriers to sos-provability that rely on the structure of a graph.

Taken as a whole, these results give some ideas of when and how to apply sos methods in attempting to prove Sidorenko's conjecture for different families of graphs.  We hope that those interested in Sidorenko's conjecture will be able to use these ideas in deciding which proof techniques are promising.  Furthermore, as each new Sidorenko graph family is discovered, its properties can again be analyzed as we have done here to determine whether or not an sos proof is at all possible, as well as where to begin such a proof.

\bibliographystyle{alpha}
\bibliography{references}

\begin{thebibliography}{HKLL16}

\bibitem[BPT12]{blekherman2012semidefinite}
Grigoriy Blekherman, Pablo~A. Parrilo, and Rekha~R. Thomas.
\newblock {\em Semidefinite {O}ptimization and {C}onvex {A}lgebraic
  {G}eometry}.
\newblock SIAM, 2012.

\bibitem[BR65]{blakleyroy}
George~R. Blakley and Prabir Roy.
\newblock H\"older type inequality for symmetric matrics with nonnegative
  entries.
\newblock {\em Proceedings of the American Mathematical Society},
  16:1244--1245, 1965.

\bibitem[BRST20]{BRST}
Grigoriy Blekherman, Annie Raymond, Mohit Singh, and Rekha~R. Thomas.
\newblock Simple graph density inequalities with no sum of squares proofs.
\newblock {\em Combinatorica}, 40:455--471, 2020.

\bibitem[CKLL18]{conlonkimleelee}
David Conlon, Jeong~Han Kim, Choongbum Lee, and Joonkyung Lee.
\newblock Some advances on {S}idorenko's conjecture.
\newblock {\em Journal of the London Mathematical Society}, 98:593--608, 2018.

\bibitem[CL21]{conlon2018sidorenko}
David Conlon and Joonkyung Lee.
\newblock Sidorenko's conjecture for blow-ups.
\newblock {\em Discrete Analysis}, 2:13 pp., 2021.

\bibitem[GL15]{gowersblog}
Tim Gowers and Jason Long.
\newblock Entropy and {S}idorenko's conjecture after {S}zegedy.
\newblock 2015.
\newblock
  \url{https://gowers.wordpress.com/2015/11/18/entropy-and-sidorenkos-conjecture-after-szegedy/}.

\bibitem[HKLL16]{KLL}
Jeong Han~Kim, Choongbum Lee, and Joonkyung Lee.
\newblock Two approaches to {S}idorenko's conjecture.
\newblock {\em Transactions of the American Mathematical Society},
  368(7):5057--5074, 2016.

\bibitem[Lon66]{london}
David London.
\newblock Inequalities in quadratic forms.
\newblock {\em Duke Mathematical Journal}, 83:511--522, 1966.

\bibitem[Lov12]{LovaszBook}
L\'{a}szl\'{o} Lov\'{a}sz.
\newblock {\em Large {N}etworks and {G}raph {L}imits}, volume~60 of {\em
  American Mathematical Society Colloquium Publications}.
\newblock American Mathematical Society, Providence, RI, 2012.

\bibitem[LS22]{li2011logarithimic}
J.L.~Xiang Li and Bal\'{a}zs Szegedy.
\newblock On the logarithimic calculus and {S}idorenko's conjecture.
\newblock {\em To appear in Combinatorica}, 2022.

\bibitem[MS59]{mulhollandsmith}
H.P. Mulholland and Cedric~A.B. Smith.
\newblock An inequality arising in genetical theory.
\newblock {\em American Mathematical Monthly}, 66:673--683, 1959.

\bibitem[Sid93]{Sid93}
Alexander~F. Sidorenko.
\newblock A correlation inequality for bipartite graphs.
\newblock {\em Graphs and Combinatorics}, 9:201--204, 1993.

\bibitem[Sze14]{sze}
Bal\'{a}zs Szegedy.
\newblock An information theoretic approach to {S}idorenko's conjecture.
\newblock {\em arXiv:1406.6738}, 2014.

\end{thebibliography}

\newpage
\section{Appendix}\label{sec:appendix}

See the last two paragraphs of Section 3 on how to read the following tables. For example, the top left sos certificate is

$$\resizebox{1cm}{!}{

} \bigg)^2\bigg]\bigg]$$

\bigskip

\begin{center}
\textbf{Sums of squares certificates for five edges}
 \end{center}

$$ %
$$

\newpage

\begin{center}
\textbf{Sums of squares certificates for six edges}
\end{center}

$$ %
$$

\newpage

\begin{center}
\textbf{Sums of squares certificates for seven edges}
\end{center}

$$ %
$$

\end{document}